\documentclass[12pt]{article}
\usepackage{e-jc}

\usepackage{amsfonts, amsmath, amssymb, amsthm}
\usepackage{graphicx}
\usepackage[caption=false]{subfig}

\usepackage[colorlinks=true,citecolor=black,linkcolor=black,urlcolor=blue]{hyperref}

\newcommand{\doi}[1]{\href{http://dx.doi.org/#1}{\texttt{doi:#1}}}
\newcommand{\arxiv}[1]{\href{http://arxiv.org/abs/#1}{\texttt{arXiv:#1}}}

\usepackage{cancel}
\usepackage{epstopdf}
\newcommand{\be}{\begin{equation}}
\newcommand{\ee}{\end{equation}}
\newcommand{\bea}{\begin{eqnarray}}
\newcommand{\eea}{\end{eqnarray}}

\newcommand{\cA}{\mathcal{A}}
\newcommand{\cM}{\mathcal{M}}

\newcommand{\cB}{\mathcal{B}}

\newcommand{\cC}{\mathcal{C}}

\newcommand{\cG}{\mathcal{G}}

\newcommand{\cI}{\mathcal{I}}

\theoremstyle{plain}
\newtheorem{theorem}{Theorem}
\newtheorem{lemma}[theorem]{Lemma}

\newtheorem{proposition}[theorem]{Proposition}

\theoremstyle{definition}
\newtheorem{definition}[theorem]{Definition}

\theoremstyle{remark}

\title{\Large Counting gluings of octahedra}

\author{Valentin Bonzom\\
\small LIPN, UMR CNRS 7030, Institut Galil\'ee, Universit\'e Paris 13, Sorbonne Paris Cit\'e\\[-0.8ex]
\small 99, avenue Jean-Baptiste Cl\'ement, 93430 Villetaneuse, France, EU\\[-0.8ex]
\small\tt bonzom@lipn.univ-paris13.fr
\and
Luca Lionni\\
\small Laboratoire de Physique Th\'eorique, CNRS UMR 8627\\[-0.8ex]
\small Universit\'e Paris XI, 91405 Orsay Cedex, France, EU and\\[-0.8ex]
\small LIPN, UMR CNRS 7030, Institut Galil\'ee, Universit\'e Paris 13, Sorbonne Paris Cit\'e\\[-0.8ex]
\small 99, avenue Jean-Baptiste Cl\'ement, 93430 Villetaneuse, France, EU\\[-0.8ex]
\small\tt luca.lionni@th.u-psud.fr
}

\begin{document}
\maketitle

\begin{abstract}
Three--dimensional colored triangulations are gluings of tetrahedra whose triangles carry the colors 0, 1, 2, 3 and in which the attaching maps between tetrahedra are defined using the colors. This framework makes it possible to generalize the notion of two--dimensional $2p$--angulations to three dimensions in a way which is suitable for combinatorics and enumeration. In particular, universality classes of three--dimensional triangulations can be investigated within this framework. Here we study colored triangulations obtained by gluing octahedra. Those which maximize the number of edges at fixed number of octahedra are fully characterized and are shown to have the topology of the 3--sphere. They are further shown to be in bijection with a family of plane trees. The enumeration is performed both directly and using this bijection. This is the first combinatorial analysis of three--dimensional colored triangulations made of building blocks which are non-melonic and have spherical topology.
\end{abstract}

\section{Introduction}

Combinatorial maps are proper embeddings of graphs into surfaces and can be thought of as discretizations of surfaces. They have been the subject of numerous studies since the seminal works of Tutte, with a special focus on the enumeration of maps. Here we are interested in the extension of this framework to three dimensions. Combinatorial maps can be built by gluing polygons: triangles, quadrangles and so on. In higher dimensions, one is interested in gluing tetrahedra or other three--dimensional building blocks.

Several approaches exist to enumerate and study properties of maps: Tutte's recursive approach \cite{CatalyticVariables}, matrix models developed by physicists \cite{MatrixReview}, bijections with labeled trees \cite{Schaeffer, Mobiles}, the topological recursion \cite{TR}, etc. Most of them do not exist or are very little developed in the three--dimensional case. For instance, the approach via matrix models generalize to tensor models \cite{GurauRyanReview, SigmaReview}. However, some powerful methods used to study matrix models, such as orthogonal polynomials, are not available to tensor models. This calls for more direct combinatorial approaches.

Interest in large random tensors (a generalization of the theory of random matrices, initially introduced in \cite{Tensors}) and enumeration of higher--dimensional ``maps'' was revived through the idea of considering colored triangulations (by triangulations in three--dimensions we mean gluings of tetrahedra). The theory of random tensors and colored triangulations has been extensively developped in the past 6 years (see the book \cite{GurauBook} for instance). Interestingly, colored triangulations had been introduced and studied in topology as a graph--theoretical representation of PL--manifolds for quite some time already (see \cite{ItalianSurvey, LinsMandel} and further references therein) under the name of graph--encoded manifolds and crystallization theory. It was further brought to our attention\footnote{We thank an anonymous referee for pointing up this part of the literature.} that there is a very close notion in combinatorics, due to Stanley, known as balanced simplicial complexes (see \cite{Balanced} for some references, a brief comparison with our framework will be offered in Section \ref{sec:Octahedra}).

We will briefly recall in Section \ref{sec:Octahedra} the equivalence between colored triangulations of PL pseudo-manifolds and a family of edge--colored graphs. This equivalence is the reason why colored triangulations seem better suited to the combinatorial analysis of three--dimensional structures than other families of gluings of tetrahedra. However, they had not been considered as combinatorial objects to be enumerated yet.

Numerous families of combinatorial maps have been defined, for instance with constraints on the degrees of faces\footnote{The word ``face'' is used in various communities to describe different objects which turn out to be relevant for our purpose (like a subsimplex on the boundary of a simplex). Here we will reserve its use to faces of combinatorial maps since its other uses admit simple way-around (e.g. triangle for a face of a tetrahedron), which is not the case for combinatorial maps.} (faces are the complements of the 1-skeleton graph and are homeomorphic to discs). Enumeration projects aim in particular at identifying the universality class of each family of maps of interest, see \cite{MatrixReview} for numerous calculations of critical exponents.

Colored triangulations provide a framework suitable to the investigation of universality classes in higher dimensions. This is because they admit natural families defined by constraints, which generalize $2p$-angulations (i.e. gluings of $2p$-gons), and the universality classes of those families can then in principle be compared to one another. In three dimensions, which is the case of interest here, colored tetrahedra can indeed be glued together to form objects with boundary which can in turn be used as new building blocks, as for instance an octahedron made of eight internal tetrahedra. Closed pseudo-manifolds are then generated by gluing copies of that octahedron. This way, families of colored triangulations built from gluings of only certain building blocks can be defined and studied. We call such building blocks \emph{bubbles} \cite{Uncoloring}.

Before universality classes can be found, it is necessary to classify colored triangulations in a similar way to the genus classification of combinatorial maps. For $p$--angulations, i.e. maps whose faces have degree $p$, the number of edges and the number of faces are not independent. Therefore, at fixed number of faces, the genus is linearly equivalent to the number of vertices of the map. In particular, planar $p$--angulations can be seen as those which maximize the number of vertices at fixed number of edges. There also is a geometrical aspect to it. Consider indeed a triangulation as made of equilateral and flat triangles. This endows the map with a geometry whose curvature lies on the vertices. In higher dimensions, following the idea of Regge calculus \cite{Regge}, a gluing of tetrahedra gets a geometry by declaring all tetrahedra regular and flat, and curvature lies around the edges. We will thus be interested in three--dimensional colored triangulations, obtained from gluing certain building blocks (octahedra, here), which maximize the number of edges.

This problem has been solved for very few families of building blocks. To grasp the difficulty, it must be noted that there \emph{is} a classification of colored triangulations with respect to Gurau's degree $\omega(T)$, by Gurau and Schaeffer \cite{GurauSchaeffer},
\begin{equation} \label{GurauDegree}
\omega(T) = \frac{3}{2} t(T) + 3 - e(T) \geq 0,
\end{equation}
where $t(T), e(T)$ are the numbers of tetrahedra and edges of the colored triangulation $T$. The main theorem \cite{1/NExpansion} is that $\omega(T)$ is a non--negative integer, providing a linear bound on the number of edges of $T$ at fixed number of tetrahedra. However, $\omega(T) = 0$ imposes severe restrictions, described at length in \cite{Melons}, which enforce a series--parallel structure on those triangulations called melonic. Moreover, melonic triangulations can only be built from specific building blocks, called melonic bubbles, which are 3--balls made of tetrahedra glued in a series--parallel way \cite{Uncoloring}. Note that a similar classification due to Fusy and Tanasa holds for a more general family of triangulations called multiorientable triangulations \cite{MOClassification}.

However, if melonic building blocks are not used, it is expected from \cite{GurauSchaeffer} that not only there are no triangulations of vanishing degree, but there are no infinite families of triangulations at fixed degree. This is because gluing non--melonic bubbles, such as octahedra, does not grow as many edges as it does with melonic ones. This is the reason why Gurau's degree is not an interesting classifying parameter for the enumeration of colored triangulations and it does not permit the study of universality classes. One instead expects the existence of a \emph{modified degree with a balance between the numbers of tetrahedra and edges which depends on the choice of bubble}. 

Such modified degrees, which admit an infinite number of triangulations for some fixed values, are particularly difficult to find in practice. So far, the existing results depend on the dimension. In even dimensions, starting from 4, there are building blocks made of only four 4-simplices for which interesting results have been obtained \cite{MelonoPlanar}. In particular, for triangulations made of these bubbles together with melonic bubbles, a general modified degree can be defined which is always non--negative. Triangulations which maximize the number of subsimples of codimension 2 are those of vanishing modified degree and their multivariate generating function has different singularities which correspond to the universality class of either random trees (when melonic bubbles dominate the asymptotics), or random planar maps (when the other bubbles dominate), or a mixed universality class (which roughly described corresponds to infinite random trees whose vertices are infinite radom planar maps). This contrast with the two--dimensional case where infinite planar $p$-angulations all lie in the same universality class.

In four dimensions, other bubbles have been investigated but only their self-gluings, i.e. the pairwise gluings of boundary tetrahedra (those triangulations are thus made of a single building block and there is no notion of universality classes). This is a generalization of unicellular maps to four dimensions which has been connected to the enumeration of meanders \cite{Meanders}.

In three dimensions, almost nothing is known beyond melonic bubbles and in fact \emph{the present article is the first identification of a modified degree for a non-melonic bubble with spherical topology}. Indeed, all bubbles made of four tetrahedra fall in the melonic class except a non-bipartite one whose boundary surface is the projective plane. Its modified degree was recently found \cite{O(N)Model}. Moreover, bipartite bubbles with six tetrahedra are all melonic except one which is a torus (the dual 1-skeleton is $K_{3,3}$). In the latter case, the modified degree was found in \cite{StuffedWalshMaps}. Therefore, in order to find non-melonic bipartite bubbles of spherical topology, one has to move up to eight tetrahedra per bubble, such as the octahedron.

Several strategies to prove the existence of modified degrees and mostly find examples where the modified degree can be written explicitly are exposed in \cite{SigmaReview}. The most promising strategy, used in four dimensions and in the 6-tetrahedra bubble mentioned above, is a bijection which turns colored triangulations into edge--colored maps with prescribed submaps (reminiscent of hypermaps). It was introduced in \cite{StuffedWalshMaps}.

An octahedron can be built by gluing eight colored tetrahedra, as shown in Figure \ref{fig:Octahedron}. It can then be used as a building block (bubble) and we will here study the triangulations obtained by gluing octahedra together. More precisely we will characterize the gluings which maximize the number of edges and deduce from them the existence of a modified degree with coefficient $5/8$ instead of $3/2$ in \eqref{GurauDegree}. We will show that they have the topology of the 3--sphere and enumerate them. As it turns out, they form a family in bijection with trees. At the end of the day, the universality class of trees is the only one which has been found in three--dimensional colored triangulations yet.

The organization is as follows. In Section \ref{sec:Octahedra} we provide all the necessary definitions. We define colored triangulations, bubbles and more particularly the octahedral bubble. We further recall that colored triangulations can be represented as regular edge--colored graphs. In Section \ref{sec:Bijection} we introduce our working representation, based on \cite{StuffedWalshMaps}, which turns colored triangulations into hybrid objects made of two types of vertices: black vertices equipped with a cyclic order of their incident edges, and square--vertices which are specific cycles of length four. The bijection is given by Theorem \ref{thm:Bijection} in Section \ref{sec:BijectionTheorem}. This representation makes it easier to count the edges of the colored triangulations. In Section \ref{sec:Maximizing} we prove our main result, Theorem \ref{thm:Dominant}, which characterizes the triangulations maximizing the number of edges, in the representation of Section \ref{sec:Bijection} and then as edge--colored graphs. Section \ref{sec:Topology} shows that they are 3--spheres. Finally, the enumeration is performed in Section \ref{sec:Enumeration}, where a simple bijection with trees is also briefly presented.

\section{Gluings of octahedra} \label{sec:Octahedra}

\subsection{Colored triangulations}

\begin{figure}
\centering
\includegraphics[scale=1]{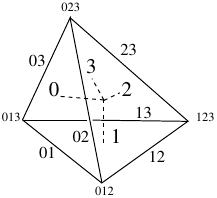}
\caption{\label{fig:ColoredTet}This is a colored tetrahedron with boundary triangles colored $0, 1, 2, 3$. Edges and vertices are respectively identified by induced pairs and triples of colors.}
\end{figure}

Colored triangulations in three dimensions consist of gluings of colored tetrahedra \cite{GurauBook, ItalianSurvey, LinsMandel, SigmaReview}. A tetrahedron is said to be colored when its triangles are labeled with the colors $0, 1, 2, 3$. Those colors are used to define the attaching maps between tetrahedra. Indeed, first observe that all sub--simplices of a single tetrahedron can be uniquely identified by sets of colors. Obviously, the triangles carry the colors 0, 1, 2, 3. Moreover, an edge is shared by two triangles, say of colors $a$ and $b$, so that the edge is identified by the pair of colors $(ab)$ (the order does not matter). Similarly, a vertex of the tetrahedron is shared by three triangles say of colors $a, b, c$ (or equivalently by three edges with pairs of colors $(ab), (bc), (ca)$) and is thus completely identified by the triple of colors $(abc)$ (again, the order does not matter). This is represented in Figure \ref{fig:ColoredTet}.

\begin{figure}
\centering
\includegraphics[scale=1]{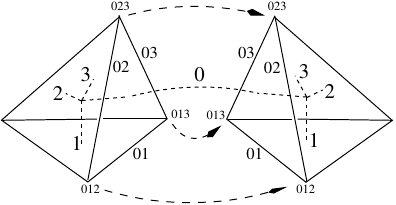}
\caption{\label{fig:GluingTet}Colors provide an unambiguous attaching constraint: one glues two tetrahedra along a triangle such that all induced colorings of the sub--simplices are preserved.}
\end{figure}

Two tetrahedra can then be glued unambiguously along a given triangle: it is the gluing which preserves all induced colorings of sub--simplices. One matches the two triangles of color $a\in\{0, 1, 2, 3\}$ of both tetrahedra so that the edges of colors $(ab)$ for all $b\neq a$ are two by two identified, and so are the three vertices of colors $(abc)$ for all $b, c \neq a$. This is pictured in Figure \ref{fig:GluingTet}.

We further restrict attention to bipartite gluings, so that tetrahedra are colored black and white and a black (white) tetrahedron is only glued to white (black) tetrahedra. Topologically, colored triangulations are homeomorphic to PL pseudo--manifolds and bipartiteness further imposes orientability of the pseudo--manifolds.

This setting, borrowed from crystallization theory and graph-encoded-manifolds introduced in topology, is also closely related to balanced simplicial complexes which have the same coloring rules. However, our colored triangulations here are typically not simplicial complexes because two tetrahedra can glued along more than a single triangle. A related aspect is that balanced simplicial complexes are typically studied with respect to the number of vertices. In our case we will study colored triangulations with respect to the number of tetrahedra (equivalently building blocks, see below) while the number of vertices is not fixed.

\begin{figure} \centering
\includegraphics[scale=.8]{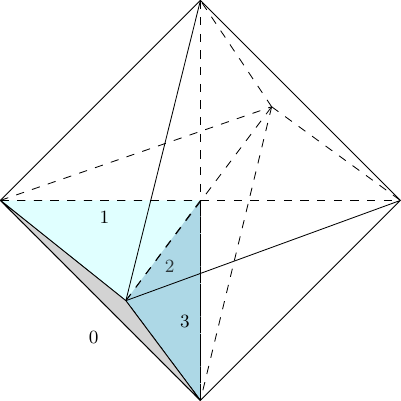}
\caption{\label{fig:Octahedron} An octahedron can be obtained by gluing eight colored tetrahedra leaving only triangles of color 0 on the boundary. Notice that the coloring is completely determined by the colors of the triangles of a single tetrahedron. Indeed, one finds the colors of the other triangles by requiring that each edge is labeled by a pair of colors.}
\end{figure}

To study the universality classes of triangulations, it is necessary to introduce subsets of triangulations. The subset of interest in the present article is based on gluings of octahedra. An octahedron can be formed by gluing eight colored tetrahedra as in Figure \ref{fig:Octahedron}, using the attaching maps described above, and such that the eight triangles on the boundary all have the color $0$. Each colored tetrahedron is thus glued to three others, along its triangles of colors 1, 2, 3. There are six inner edges. Following Figure \ref{fig:Octahedron}, the octahedron can be thought of as a gluing of two up and down pyramids. The base of the pyramids consists of triangles of color 1 with four inner edges of colors $(12)$ and $(13)$ alternating. The heights of both pyramids are edges of colors $(23)$. The octahedron is symmetric under permutations of the colors 1, 2, 3.

\begin{definition}
The set $\cG$ denotes the subset of connected colored triangulations obtained by gluing copies of the octahedron of Figure \ref{fig:Octahedron}.
\end{definition}


\subsection{Colored graphs}

Colored triangulations have a major combinatorial advantage: they can be represented by graphs with colored edges. The graph corresponding to a triangulation is simply obtained by taking the 1--skeleton of the dual. The 1--skeleton of the dual to a colored tetrahedron is simply a vertex with four incident edges (one for each triangle). Each edge is dual to a triangle and therefore comes with a color, as represented by the dashed edges of Figure \ref{fig:ColoredTet}. As shown in Figure \ref{fig:GluingTet}, the gluing of two tetrahedra along a triangle of color $a$ is represented dually by an edge joining the two vertices dual to the tetrahedra. That edge carries the color $a$ of the triangle.

This way one obtains graphs which are: bipartite with black and white vertices, regular of degree $4$ (since a tetrahedron has four triangles), and such that the four colors are incident exactly once to each vertex. We simply call those graphs \emph{colored graphs} in the remaining. Thanks to the edge coloring, a colored graph has all the information necessary to reconstruct the corresponding colored triangulation, called the pseudo--complex associated to the graph in crystallization theory. In particular, cells of dimension $3-k$ labeled by $k$ colors $a_1, \dotsc, a_k$ are represented in the graph as maximally connected components whose edges only have the colors $a_1, \dotsc, a_k$ (called $k$-residues in crystallization theory).

A \emph{bubble} is a maximally connected sub--graph with the colors 1, 2, 3, but not 0. Obviously, all colored graphs are obtained by considering arbitrary (finite) subsets of bubbles and adding the color 0 to connect white to black vertices. Colored graphs with constrained bubbles are a generalization of combinatorial maps with constrained face degrees. A bubble can be seen as a building block and one is interested in the universality classes of the triangulations obtained by gluing bubbles which are chosen in a finite set.

The octahedron of Figure \ref{fig:Octahedron} can be described as a bubble denoted $B$, by taking the 1--skeleton of the dual and ignoring the triangles of color $0$ (since they are not glued), 
\begin{equation}
\begin{array}{c} \includegraphics[scale=.75]{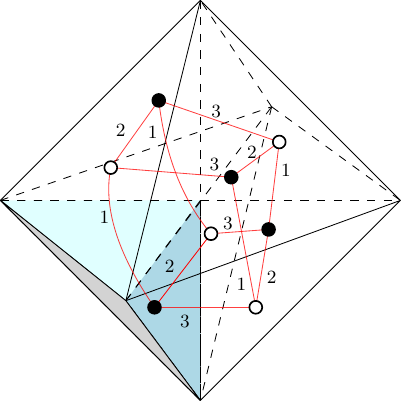} \end{array} 
\qquad \Rightarrow \qquad B = \begin{array}{c} \includegraphics[scale=.75]{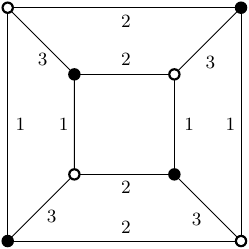} \end{array}
\end{equation}
The bubble $B$ simply comes from redrawing the red edges in a more readable fashion.

Bubbles are the natural building blocks of colored triangulations\footnote{They are in a one-to-one correspondence with $U(N)^3$-invariant polynomials in the entries of a complex tensor $T_{a_1 a_2 a_3}$ \cite{Uncoloring}. This enables the generalization of the famous relationship between random matrices and combinatorial maps to random tensors and colored triangulations.}. Let us look at the simplest bubbles before we focus on the octahedral one. With two vertices, there is a single bubble which corresponds to the gluing of two tetrahedra along their three triangles of colors 1, 2, 3. In the introduction we insisted that very little is known in three dimensions beyond melonic bubbles. A melonic bubble can be defined as a bubble for which there is a sequence of contractions of parallel edges which reduces it to the two-vertex bubble. It is easy to see \cite{Uncoloring} that all (bipartite) bubbles on four vertices and six vertices are melonic expect one whose graph is $K_{3,3}$, is topologically a torus and was studied in \cite{StuffedWalshMaps}. Therefore if one wants to study a non-melonic bubble with spherical topology, it is necessary to look at bubbles with eight vertices and the octahedron is a natural example.

\begin{proposition}
The set $\cG$ of gluings of octahedra can be equivalently thought of as the set of connected colored graphs built from gluing copies of $B$ via edges of color 0 between all black and white vertices\footnote{The bubbles and vertices of each graph are not labeled.}.
\end{proposition}

\subsection{Bicolored cycles and enhanced degree}


From now on we call colored triangulations or graphs the elements of $\cG$. Our goal is to find those which maximize the number of edges as triangulations, and for fixed number of octahedra. Edges of a triangulation have pairs of colors. Each octahedron brings six edges with labels $(ab)$ for $a, b \in\{1, 2, 3\}$. The remaining edges have color type $(0a)$ for $a=1, 2, 3$. They are the ones which are not fixed by the number of octahedra and that we want to maximize.

By duality, edges of the triangulation are represented by maximally connected components with pairs of colors in the colored graphs. They are obviously bicolored cycles with alternating colors. Therefore, we want to maximize the number of bicolored cycles with colors $(0a)$ at fixed number of bubbles.

\begin{definition} \label{def:Faces}
For $G\in\cG$ we denote $F_a(G)$ the number of bicolored cycles of colors $(0a)$, $a\in\{1, 2, 3\}$, and $F(G) = F_1(G) + F_2(G) + F_3(G)$ the total number of such bicolored cycles. Maximizing the number of edges of a triangulation is equivalent to maximizing $F(G)$ of the corresponding colored graph.
\end{definition}

Let us look at the colored graphs with a single bubble. They are obtained by adding edges of color $0$ between the black and white vertices of $B$. There are $4! = 24$ ways to do so (ignoring symmetries), but not all of them maximize $F(G)$. It can be checked explicitly that in order to maximize the number of bicolored cycles, one must add the edges of color 0 parallel to the edges of a chosen color $a = 1, 2, 3$. Consequently, three graphs with a single bubble maximize the number of bicolored cycles,
\begin{equation} \label{CubePairing}
\begin{array}{c} \includegraphics[scale=.8]{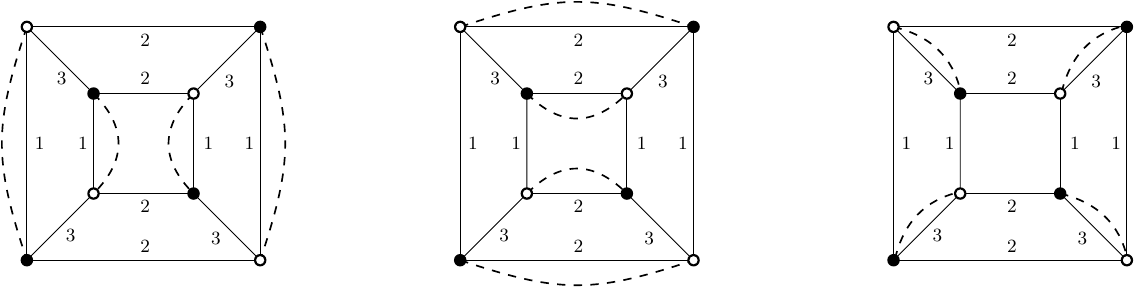} \end{array}
\end{equation}
From here onwards, edges of color 0 are represented with dashed edges. If the edges of color 0 are parallel to those of color $a$, there are 4 bicolored cycles of colors $(0a)$, 2 of colors $(0b)$ and 2 of colors $(0c)$ for the complementary colors $\{b,c\} = \{1, 2, 3\} \setminus \{a\}$, leading to a total of 8.

For colored triangulations $T$ built from gluings of a given bubble $B$ made of $2p$ tetrahedra, Gurau's degree \eqref{GurauDegree} can be adapted and takes the form
\begin{equation} \label{GurauDegreeBubble}
\omega_B(T) = \frac{p-1}{p} t(T) + 3 - e_{0}(T),
\end{equation}
where $e_0(T)$ is the number of edges of color type $(0a)$ for all $a\in\{1, 2, 3\}$ (i.e. the number of bicolored cycles of the corresponding colored graph). However, when $B$ is non--melonic, $T$ cannot grow as many edges as in a melonic triangulation. One therefore expects the existence of $s_B < (p-1)/p$ such that the modified degree
\begin{equation} \label{OctahedraDegree}
\tilde{\omega}_B(T) = s_B\,t(T) + 3 - e_0(T)
\end{equation}
is still a non--negative integer. With eight tetrahedra, $p=4$ and $(p-1)/p = 3/4$, while we will find that $\tilde{\omega}_B(T)$ is still non--negative for $s_B = 5/8 < 3/4$. In particular for the graphs in \eqref{CubePairing}, the number of tetrahedra is $t(T) = 8$ and $e_0(T) = F(T) = 8$, confirming $\tilde{\omega}_B(T) = 0$. Moreover, all $s > 5/8$ would make $\tilde{\omega}(T)$ negative for an infinite number of $T$. This will be a consequence of Theorem \ref{thm:Dominant}.

\section{Another representation} \label{sec:Bijection}

\subsection{The new objects} 

A map is a graph equipped with a rotation system, i.e. a cyclic ordering of the edges incident to each vertex. A corner is a pair of consecutive edges around a vertex and it can be oriented say counter--clockwise. This gives rise to the notion of \emph{faces}, which are cycles obtained by following the edges and the corners to go from one edge to another at a vertex. Equivalently, a map can be thought of as a properly embedded graph (without crossings) on a surface, up to isotopy, and the faces are then the connected components of complement of the graph. Each face is homeomorphic to a disc.

\begin{definition} \label{def:M}
We consider the following set $\cM$ of connected graphs made of
\begin{itemize}
\item cycles of length four with alternating colors 1, 2, 1, 2. We call those cycles square--vertices.
\item black vertices of arbitrary (finite) degrees, every one of them being equipped with a cyclic ordering of its incident edges,
\item edges connecting black vertices to square--vertices, so that the four vertices inside each square--vertex are connected to black vertices.
\end{itemize}
We call the edges and vertices of the square--vertices inner edges and inner vertices and we simply use edges for those which connect black vertices to inner vertices.
\end{definition}

While the orders of the incident edges at black vertices matter, there is no order around inner vertices. We can therefore always represent square--vertices as squares with four exterior edges.
\begin{equation} \label{SquareVertex}
\begin{array}{c} \includegraphics[scale=.7]{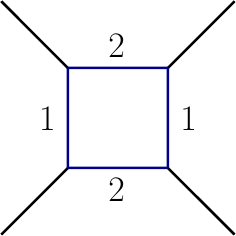} \end{array}
\end{equation}

An element $M\in\cM$ is a sort of hybrid between a graph and a map. It has two distinguished subgraphs which are maps and come from the colored structure of square--vertices. The map $M^{(1)}$ of color 1 (respectively $M^{(2)}$ of color 2) is obtained by deleting all inner edges of color 2 (respectively of color 1) from all square--vertices of $M$. $M^{(1)}$ and $M^{(2)}$ are maps because after deleting the edges of color 2 (or 1), each inner vertex is bivalent and there is a unique cyclic order between its incident edge and inner edge. 


\begin{definition} \label{def:FacesM}
The faces of $M\in\cM$ are defined as the faces of the combinatorial maps $M^{(1)}$ and $M^{(2)}$ as well as the black vertices.
\end{definition}

The map $M^{(1)}$ (or $M^{(2)}$) consists of the edges of $M$ and the inner edges of color 1 (or 2). We can thus think of the edges of $M$ as carrying the pair of colors $(12)$ and of $M^{(1)}$ (or $M^{(2)}$) as the map which has all the edges carrying the color 1 (or 2). The reason we consider black vertices as faces in addition is that no edges nor inner edges carry the color 3. Thus we could define a map $M^{(3)}$ obtained by erasing all edges and inner edges which do not have the color 3, i.e. all of them. This way $M^{(3)}$ would consist of isolated black vertices, each with an external face around it.

As outlined before, after deleting the inner edges of a given color, the inner vertices become bivalent. We can therefore remove them and merge their incident edge and inner edge. This way, the square--vertex of \eqref{SquareVertex} becomes $\begin{array}{c} \includegraphics[scale=.3]{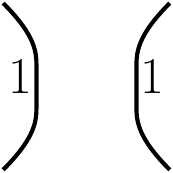} \end{array}$ in $M^{(1)}$ and $\begin{array}{c} \includegraphics[scale=.3]{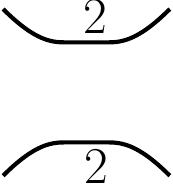} \end{array}$ in $M^{(2)}$. As a consequence, the map $M^{(c)}$, $c=1, 2$, has $2b$ edges and $V$ vertices, if $b$ is the number of square--vertices of $M$ and $V$ its number of black vertices.\\

Let us give some examples with a single square--vertex. Those which maximize the number of faces are the following
\begin{equation} \label{1SquareVertexEx}
\begin{array}{c} \includegraphics[scale=.75]{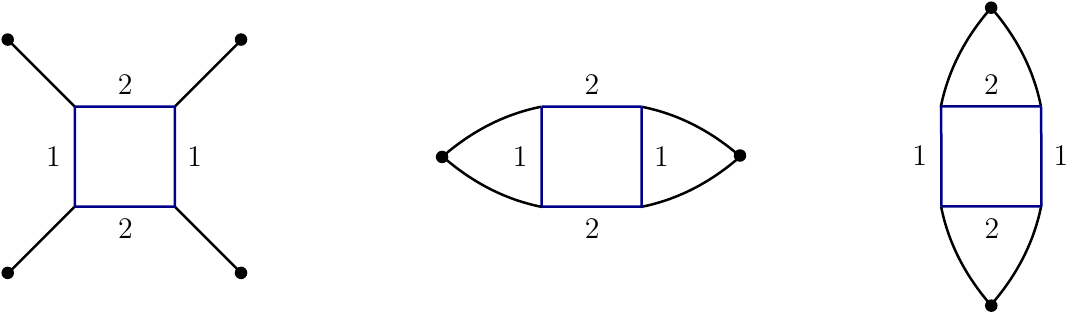} \end{array}
\end{equation}
The left one is a tree (with one square--vertex and four leaves). Its submap $M^{(1)}$ has two disconnected edges of color 1, hence two faces of color 1, and similarly for the color 2. There are four black vertices, hence a total of eight faces. The middle one has four faces of color 1 ($M^{(1)}$ has two disconnected loops), two faces of color 2 ($M^{(2)}$ has two vertices and two parallel edges joining them), plus two black vertices, giving eight faces again. The one on the right is obtained by permuting the colors 1 and 2.

By contrast, the following element of $\cM$ only has six faces,
\begin{equation}
M = \begin{array}{c} \includegraphics[scale=.65]{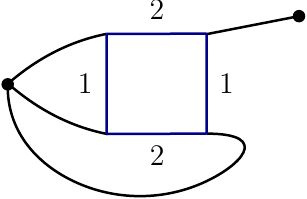} \end{array} \quad M^{(1)} = \begin{array}{c} \includegraphics[scale=.65]{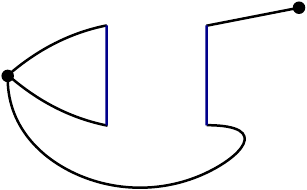} \end{array} \quad M^{(2)} = \begin{array}{c} \includegraphics[scale=.65]{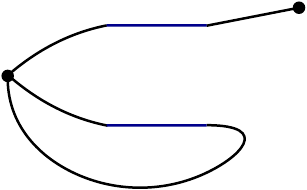} \end{array}
\end{equation}
since both $M^{(1)}$ and $M^{(2)}$ consist of an edge and a loop attached at one end, giving rise to two faces (plus two black vertices).

\subsection{The bijection} \label{sec:BijectionTheorem}

This bijection is a particular case of a generic bijection exposed in \cite{StuffedWalshMaps}. We adapt the bijection and the proof to make the article both self-contained and simpler to read than \cite{StuffedWalshMaps}.

\begin{theorem} \label{thm:Bijection}
There is a bijection between $\cG$ and $\cM$ which maps bubbles to square--vertices, bicolored cycles of colors $(01)$ to (respectively of colors $(02)$) to faces of $M^{(1)}$ (respectively $M^{(2)}$) and bicolored cycles of colors $(03)$ to black vertices (we recall that the definition of the faces of $M$ are in \ref{def:FacesM}).
\end{theorem}

The proof of the bijection itself is particularly simple and can be reduced to the following construction.
\begin{itemize}
\item For a graph $G\in\cG$, contract the edges of color 3 and identify the resulting parallel edges of the same color. This way, each bubble reduces to a square with alternating colors 1, 2. This will be a square--vertex in $M\in\cM$.
\item Moreover, each vertex, which results from the merging of a black and white vertices, is now incident to exactly two edges of color 0 (to remember which one was incident to the white vertex, it is sufficient to orient edges of color 0 from, say, white to black). Edges of color 0 thus forms a disjoint set of unicycles. Each of them can be represented as a black vertex, where the cyclic order of the edges around a cycle now translates into a cyclic order around the black vertex. This gives rise to the black vertices of $M\in\cM$.
\end{itemize}
Going from $\cM$ to $\cG$ just follows the reverse process.

The most interesting part of the bijection still concerns the bicolored cycles and the faces. It is much more convenient to track them down using a labelling of the vertices and half--edges. While this will make the proof below seemingly more complicated, because the labelling forces us to consider the symmetries of the relabellings, we think it offer several advantages. First, we will use permutations to describe both the colored graphs and the maps. They act on the vertex set for colored graphs and on the half-edge set for maps. Those encodings are really standard in both cases and it is illuminating to see without much effort that they actually coincide. A second aspect is that again for both colored graphs and maps, those encodings can be used as data structures to generate colored graphs (see \cite{Catalogues} and much more cataloguing and programming at \url{http://cdm.unimo.it/home/matematica/casali.mariarita/DukeIII.htm}) and maps \cite{FusyThese}.

{\bf Proof.} We are going to show that elements of $\cG$ and $\cM$ can be described by the same system of permutations with the same symmetries. Let us start with $\cG$.

\paragraph{Elements of $\cG$ as orbits of permutations.}
Let $G\in\cG$ with $b$ bubbles (i.e. $b$ copies of $B$). It can be described through one permutation per color which sends white vertices to black vertices. To make it explicit, consider a labeling of the vertices of $B$, $\{1_\circ, 1_\bullet, 2_\circ, 2_\bullet, 3_\circ, 3_\bullet, 4_\circ, 4_\bullet\}$ where 
\begin{itemize}
\item $i_\circ, i_\bullet$ are respectively white and black vertices connected by the color 3,
\item the vertices $1_\circ, 1_\bullet$ are connected to $2_\bullet, 2_\circ$ by edges of color 1,
\item the vertices $1_\circ, 1_\bullet$ are connected to $4_\bullet, 4_\circ$ by edges of color 2,
\end{itemize}
\begin{equation}
\begin{array}{c} \includegraphics[scale=.75]{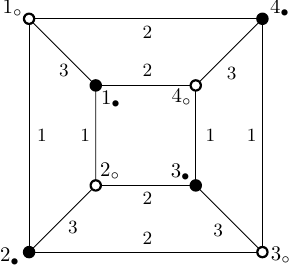} \end{array}
\end{equation}
We refer to this labeling of $B$ as the \emph{standard} labeling.

The bubble $B$ is then encoded in three permutations $\tau^{(1)}, \tau^{(2)}, \tau^{(3)}$. The permutation $\tau^{(c)}$ maps each white vertex to the black vertex it is connected to by the color $c = 1, 2, 3$. Using the cycle notation for permutations,
\begin{equation}
\tau^{(1)} = (12)(34), \qquad \tau^{(2)} = (14)(23), \qquad \tau^{(3)} = \operatorname{id}
\end{equation}

Label now the bubbles of $G$ from $1$ to $b$. Using a standard labeling on each copy of $B$ in $G$ induces a labeling of the vertices of $G$, $\{i_{\circ n}, i_{\bullet n}\}$ for $i = 1, 2, 3, 4$ and $n = 1, \dotsc, b$. We also denote $V = \{i_n\}$ the set of labels and $|V| = 4b$ is the number of, say, white vertices of $G$. The structure of $G$ is then determined by permutations $\tau^{(c)}_G : V \to V$ which are obtained by concatenating the $b$ copies of $\tau^{(c)}$ for each color $c=1, 2, 3$,
\begin{equation} \label{FixedPermutations}
\tau^{(1)}_G = \prod_{n=1}^b (1_n 2_n) (3_n 4_n), \qquad \tau^{(2)}_G = \prod_{n=1}^b (1_n 4_n) (2_n 3_n), \qquad \tau^{(3)}_G = \operatorname{id}_V 
\end{equation}
and adding a permutation $\tau^{(0)}_G:V\to V$ which maps the label of a white vertex to the label of the black vertex it is connected to by an edge of color $0$ in $G$.

Clearly, the four permutations $(\tau^{(0)}_G, \tau^{(1)}_G, \tau^{(2)}_G, \tau^{(3)}_G)$ determine $G$. However, $G$ can be described by several such sets of permutations. Indeed, for each copy of $B$ in $G$, there are four standard labelings determined by the four possible choices for, say, the first white vertex $1_{\circ n}$ of the $n$--th bubble. There is thus a group of symmetry which transforms $(\tau^{(0)}_G, \tau^{(1)}_G, \tau^{(2)}_G, \tau^{(3)}_G)$ through those re--labelings. For the bubble $n \in \{1, \dotsc, b\}$, if $1_{\circ n}$ is redefined to be $2_{\circ n}$, then the labeling being standard imposes 
\begin{equation}
1_{\circ n} \to 2_{\circ n}, \qquad 2_{\circ n} \to 1_{\circ n}, \qquad 3_{\circ n} \to 4_{\circ n}, \qquad 4_{\circ n} \to 3_{\circ n}.
\end{equation}
It is thus generated by $\gamma_{2n} = (1_n 2_n)(3_n 4_n)$ (and the identity on the labels of the other bubbles), which acts on $(\tau^{(0)}_G, \tau^{(1)}_G, \tau^{(2)}_G, \tau^{(3)}_G)$ by conjugation on the four permutations. Since the labeling of the bubble $n$ is still standard, the permutations $\tau^{(1)}_G, \tau^{(2)}_G, \tau^{(3)}_G$ are unchanged (as can be checked by direct calculation). However, the conjugation affects $\tau^{(0)}_G$. 

If one redefines $1_{\circ n}$ to be $3_{\circ n}$, the transformation is generated by $\gamma_{3n} = (1_n 3_n) (2_n 4_n)$ which also commutes with, and thus does not affect $\tau^{(1)}_G, \tau^{(2)}_G, \tau^{(3)}_G$. Similarly $\gamma_{4n} = (1_n 4_n) (2_n 3_n)$ generates the transformation to the standard labeling where $1_{\circ n}$ becomes $4_{\circ n}$ and does not change  $\tau^{(1)}_G, \tau^{(2)}_G, \tau^{(3)}_G$. The group
\begin{equation} \label{TransfoGroup}
\Gamma = \Gamma_1 \times \dotsb \times \Gamma_b \qquad \text{with $\Gamma_n = < \gamma_{2n}, \gamma_{3n}, \gamma_{4n} >$ for $n=1, \dotsc, b$}
\end{equation}
is the group of standard relabelings of the bubbles of $G$. It is a subgroup of $S_{|V|}$ the group of permutations on $V$.

In addition, the labeling of the bubbles from $1$ to $b$ can be changed at fixed $G$ while preserving the standard labelings. These transformations are generated by the symmetric group $S_b$ which acts as follows
\begin{equation}
\pi \in S_b: V \to V, \qquad i_n \mapsto i_{\pi(n)}
\end{equation}
which makes it a subgroup of $S_{|V|}$. Again, those transformations do not change $\tau^{(1)}_G$, $\tau^{(2)}_G$, $\tau^{(3)}_G$ which are therefore independent of $G$ and always given by \eqref{FixedPermutations}. They are instead characteristic of the set $\cG$ and we denote them $\tau^{(1)}_\cG, \tau^{(2)}_\cG, \tau^{(3)}_\cG$ from now on.

$S_b$ and $\Gamma$ are both subgroups of $S_{|V|}$ and their product also is a subgroup, the one generated by $S_b$ and $\Gamma$, which we denote
\begin{equation}
\widetilde{\Gamma} = < S_b, \Gamma >.
\end{equation}
A graph $G$ is thus an orbit $\{ \gamma \tau^{(0)} \gamma^{-1} \}_{\gamma\in \widetilde{\Gamma}}$, with an additional constraint on $\tau^{(0)}\in S_{|V|}$ so that it reconstructs a connected $G\in\cG$. Connectedness of $G$ means that there is a path going from white to black to white vertices and so on between any two vertices of $G$. Going from a white to a black vertex is performed by an iteration of one of the four permutations $\tau^{(0)}_G, \tau^{(c)}_\cG$, $c=1, 2, 3$ (and from a black to a white vertex is performed by the inverse permutations). Therefore the constraint on $\tau^{(0)}$ is that the group generated by $\tau^{(0)}, \tau^{(1)}_\cG, \tau^{(2)}_\cG, \tau^{(3)}_\cG$ acts transitively on $V$.

Furthermore, it is immediate to see that the bicolored cycles of colors $(0a)$ for $a=1, 2, 3$ are the cycles of $\tau^{(0)}_G \tau^{(a)-1}_\cG$. The number of such cycles $C(\tau^{(0)}_G \tau^{(a)-1}_\cG)$ is clearly invariant along the orbit generated by $\widetilde{\Gamma}$, as expected and
\begin{equation}
F(G) = C(\tau^{(0)}_G \tau^{(1)-1}_\cG) + C(\tau^{(0)}_G \tau^{(2)-1}_\cG) + C(\tau^{(0)}_G \tau^{(3)-1}_\cG).
\end{equation}

\paragraph{Elements of $\cM$ as orbits of permutations.}
Let $M\in\cM$ and label the square--vertices of $M$ from $1$ to $b$. We define a standard labeling of $M$ as a labeling of the edges of $M$ such that for each square--vertex of label $n$,
\begin{itemize}
\item the four incident edges receive the labels $(1_n, 2_n, 3_n, 4_n)$
\item $1_n$ is connected to $2_n$ by an inner edge of color 1,
\item $1_n$ is connected to $4_n$ by an inner edge of color 2.
\end{itemize}
The standard labeling puts for instance the label $1_n$ on the top--left edge of \eqref{SquareVertex} and $2_n, 3_n, 4_n$ counter--clockwise. The edges of $M$ are thus labeled by the elements of $V = \{i_n\}_{i\in\{1, 2, 3, 4\}, n=1, \dotsc, b}	$.

We define $\sigma_M$ as for ordinary combinatorial maps. Each black vertex of $M$ is equipped with a cyclic order of its incident edges. This gives a cyclic order of the corresponding labels which provides a cycle to $\sigma_M$.

While $\sigma_M$ is sufficient to reconstruct $M$, several permutations correspond to the same $M$. For each square--vertex, there are four standard labelings. From a given one, the other three are obtained by permuting the labels $(1_n, 2_n, 3_n, 4_n)$ via $\gamma_{2n}, \gamma_{3n}$ and $\gamma_{4n}$ precisely. Those permutations generate $\Gamma_n$ and act by conjugation on $\sigma_M$, as well as $\Gamma$ defined in \eqref{TransfoGroup}. One can also relabel the $b$ square--vertices using a permutation in $S_b$. Therefore the group of symmetry is $\widetilde{\Gamma}$ and $M$ is an orbit $\{\gamma \sigma_M \gamma^{-1}\}_{\gamma\in\widetilde{\Gamma}}$, with an additional constraint on $\sigma_M\in S_{|V|}$ to ensure connectedness. Let us assume for a moment that the connectedness constraint is the same as for $\tau^{(0)}$ for graphs in $\cG$. This provides a bijection between $\cG$ and $\cM$.

The permutations $\sigma_M$ and $\tau^{(0)}_G$ are defined over $V = \{i_n\}_{i\in\{1, 2, 3, 4\}, n=1, \dotsc, b}$ where $b$ is in one case the number of square--vertices and in the other case the number of bubbles. The number of square--vertices in $M$ thus corresponds to the number of bubbles in $G$.

Black vertices are in correspondence with cycles of $\sigma_M$. Since $\tau^{(3)}_\cG$ is the identity on $V$, one can interpret black vertices as cycles of $\sigma_M \tau^{(3)-1}_\cG$. By replacing $\sigma_M$ with a corresponding $\tau^{(0)}_G$ (in the same orbit and thus having the same cycle structure), one finds that the black vertices of $M$ correspond to the bicolored cycles of colors $(03)$ of $G$.

One needs to identify the faces of the submaps $M^{(1)}, M^{(2)}$. Faces of $M^{(1)}$ follow consecutively the corners of black vertices and the corners of color 1 of the square--vertices. Corners of black vertices correspond to the action of $\sigma_M$. Corners of color 1 on square--vertices are obtained by deleting the inner edges of color 2 in each square--vertex. These corners are thus given by the permutation $(1_n 2_n) (3_n 4_n)$ for the square--vertex $n$. Therefore, a face of $M^{(1)}$ is a cycle of $\sigma_M \tau^{(1)-1}_\cG$. Similarly, the faces of $M^{(2)}$ are the cycles of $\sigma_M \tau^{(2)-1}_\cG$. Replacing $\sigma_M$ with $\tau^{(0)}_G$ in the same orbit does not change the number of such cycles. This gives a one--to--one correspondence between faces of $M^{(1)}$ and bicolored cycles of colors $(01)$ in $G$ (and faces of $M^{(2)}$ and bicolored cycles of colors $(02)$ in $G$).

Finally, connectedness of $M$ means that there is a path between any two edges of $M$. Such paths consist of corners around black vertices (action of $\sigma_M$) and corners of square--vertices (action of $\tau^{(1)}_{\cG}$ and $\tau^{(2)}_{\cG}$). The orbit $\{\gamma \sigma \gamma^{-1}\}_{\gamma \in \widetilde{\Gamma}}$ thus reconstructs a connected $M$ if and only if the group generated by $\sigma, \tau^{(1)}_\cG, \tau^{(2)}_\cG$ acts transitively on $V$. Together with the fact that $\tau^{(3)}_\cG = \operatorname{id}_V$, this is a the same constraint as for $\cG$. \qed

\section{Maximizing the number of faces} \label{sec:Maximizing}

\subsection{Trees}

\begin{definition}
We say $M\in\cM$ is a tree if its edges (not inner edges) are all bridges. Equivalently, it is a tree if the graph obtained by replacing square--vertices with ordinary vertices of degree four is a tree in the ordinary sense.
\end{definition}

\begin{proposition} \label{prop:Trees}
Let $T\in\cM$ be a tree with $b(T)$ square--vertices. Its number of faces is $F(T) = 5b(T) + 3$.
\end{proposition}

{\bf Proof.} The formula is true for the unique $T$ with $b(T) = 1$, as shown in \eqref{1SquareVertexEx}. Let $T\in\cM$ be a tree with $b(T) \geq 2$ and assume that the formula is true for all $b(T')<b(T)$. $T$ being a tree contains a square--vertex $v_{\Box}$ incident to three leaves $v_1, v_2, v_3$, and a fourth edge incident to a black vertex $\bar{v}$ which is not a leaf. Consider the tree $T'$ obtained from $T$ by removing $v_\Box$ and its four incident edges. It satisfies $b(T') = b(T)-1$, so that $F(T') = 5b(T) - 2$ from the induction hypothesis. We then compare directly the number of faces lost from $T$ to $T'$. The map $T^{(1)}$ has a face going around, say, $v_1$ and $v_2$, and another one passing by $\bar{v}$ and all around $v_3$ back to $\bar{v}$, as indicated by the dotted lines below.
\begin{equation}
\begin{array}{c} \includegraphics[scale=.75]{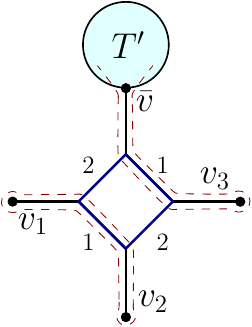} \end{array}
\end{equation}
The map $T'^{(1)}$ loses the face around $v_1$ and $v_2$ but the other one just gets shorter as it does not go around $v_3$ anymore. Therefore $F(T'^{(1)}) = F(T^{(1)}) - 1$, and similarly for $T^{(2)}$. Each black vertex being counted as a face, three of them are lost from $T$ to $T'$, so that $F(T) = F(T') + 5$, which concludes the induction. \qed

\subsection{Edge unhooking}

If an edge is not a bridge, it can be unhooked in a unique way from its incident black vertex $\bar{v}$, thereby creating a new black vertex $\bar{v}'$ of degree one, decreasing the degree of $\bar{v}$, and without increasing the number of connected components.
\begin{equation}
\begin{array}{c} \includegraphics[scale=0.4]{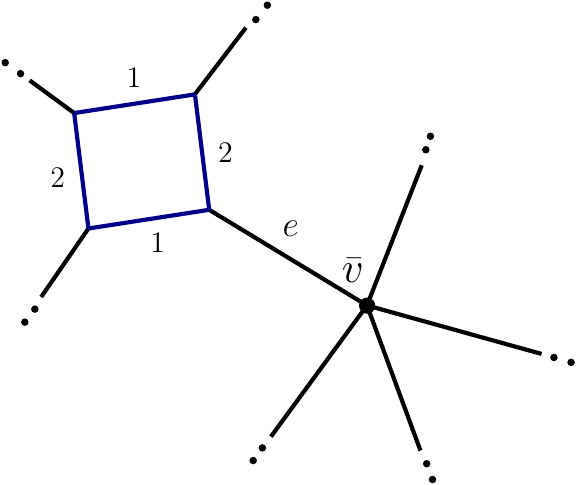} \end{array} \qquad \to \qquad \begin{array}{c} \includegraphics[scale=0.4]{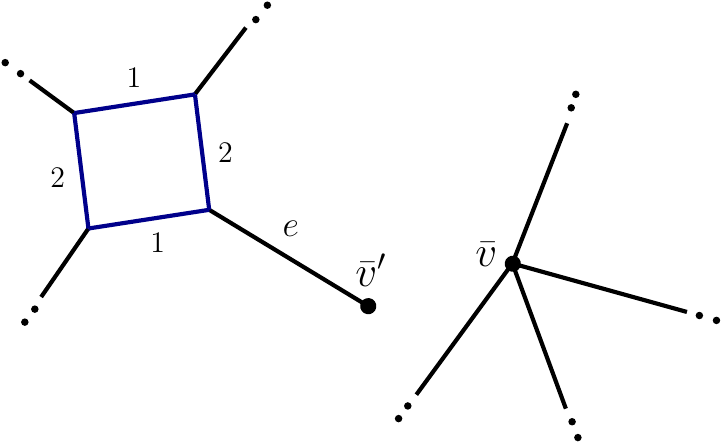} \end{array}
\end{equation}

Let us denote $\cI_2(e) \subset \{1, 2\}$ the set of colors $c$ for which the edge $e$ is incident to two faces in the map $M^{(c)}$ ($e$ can be incident to one or two faces in $M^{(c)}$), and $|\cI_2(e)|$ the number of colors in that set.

\begin{lemma} \label{lemma:Diff}
Let $M\in \cM$ and $e$ an edge of $M$ that is not a bridge, and $M_e\in\cM$ obtained by unhooking $e$. Then the variation of the number of faces is
\begin{equation}
\Delta(M, e) \equiv F(M_e) - F(M) = 3 - 2|\cI_2(e)|.
\end{equation}
\end{lemma}

{\bf Proof.} If there is a single face of color $c\in\{1, 2\}$ incident to $e$ in $M$, then this face splits into two in $M_e$. If there are two faces incident to $e$ in $M$, then they merge in $M_e$. There is moreover one additional black vertex in $M_e$. \qed

If $|\cI_2(e)| \leq 1$, then $\Delta(M, e) >0$, meaning that unhooking $e$ increases the number of faces of $M$. However, if $\cI_2(e) = \{1, 2\}$, then $\Delta(M, e) = -1$ and a face is lost upon unhooking $e$.

\subsection{The dominant gluings}

\begin{definition}
We say that $M\in\cM$ is dominant, or maximal, if it maximizes the number of faces among maps with the same number of square--vertices. We denote $\cM_{\max}$ the subset of dominant elements of $\cM$.
\end{definition}

To discuss the dominant elements of $\cM$, it will be useful to adapt the usual definition of a bond to our context, using unhooking instead of edge--deletion. A \emph{$k$--bond} is a minimal edge-cut comprised of $k$ edges, i.e. a set $S$ of edges such that unhooking all of them disconnects a connected graph into two connected components while unhooking the edges of any proper subset of $S$ does not.

\begin{figure} \centering
\subfloat[\label{fig:2Bond}Two 2--bonds whose edges are connected by the color $b$ of an inner edge.]{\includegraphics[scale=0.35]{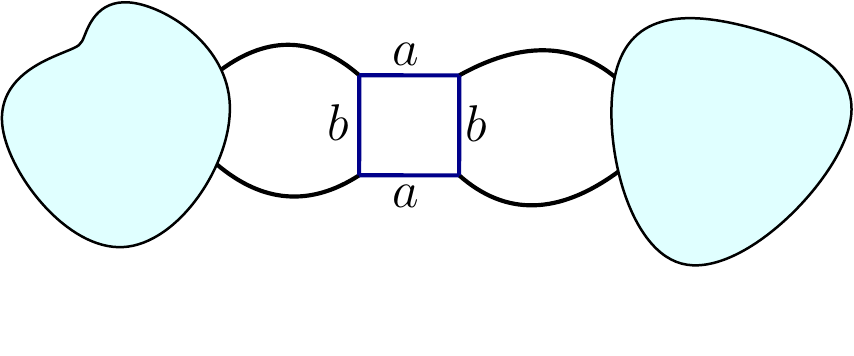}}{\qquad}
\subfloat[\label{fig:2BondBis}Two 2--bonds involving opposite edges of the square--vertex.]{\includegraphics[scale=0.35]{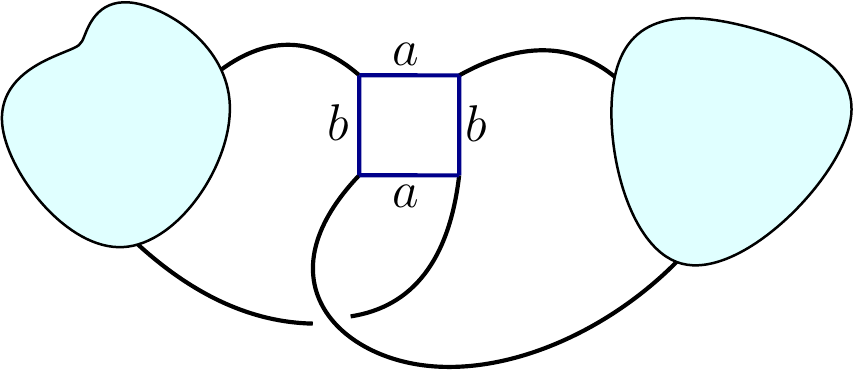}}{\qquad}
\subfloat[\label{fig:4Bond}The 4--bond case.]{\includegraphics[scale=0.35]{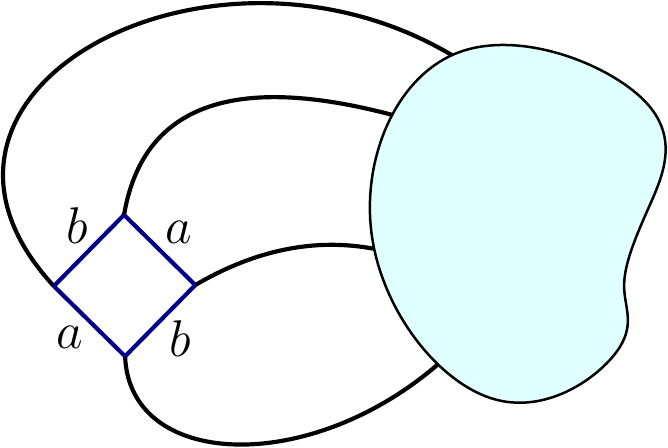}}
\caption{\label{fig:bonds}The possible bond decompositions around a square--vertex.}
\end{figure}

\begin{proposition} \label{prop:2Bond}
If $M\in\cM_{\max}$, then the four edges incident to every square--vertex of $M$ are either bridges or form two 2--bonds like in Figure \ref{fig:2Bond}. Equivalently, 2--bonds and 4--bonds like in Figure \ref{fig:2BondBis}, \ref{fig:4Bond} cannot occur.
\end{proposition}

{\bf Proof.} From the Lemma \ref{lemma:Diff} an edge $e$ in $M$ is either a bridge or satisfies $\Delta(M, e) = -1$ (or else $M$ would not be dominant). In this proof, we denote $e, e_1, e_2, e_3$ the four edges incident to a square--vertex, such that there is an inner edge of color 1 from $e$ to $e_1$ and an inner edge of color 2 from $e$ to $e_2$,
\begin{equation}
\begin{array}{c} \includegraphics[scale=.5]{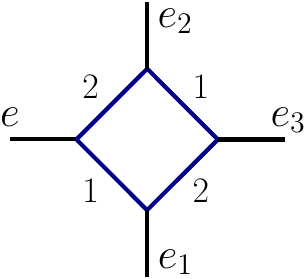} \end{array}
\end{equation}

Assume that $e$ is a bridge, 
\begin{equation}
\begin{array}{c} \includegraphics[scale=.5]{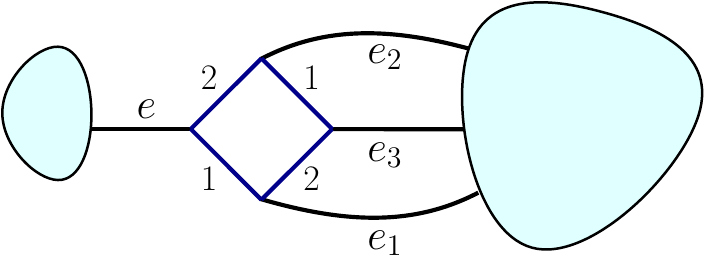} \end{array}
\end{equation}
so that $\cI_2(e) = \emptyset$, meaning that it is incident to a single face in $M^{(1)}$ and in $M^{(2)}$. This implies that there is also a single face incident to $e_1$ in $M^{(1)}$. Therefore $|\cI_2(e_1)| \leq 1$. If $e_1$ is not a bridge, then $\Delta(M, e_1)>0$ and unhooking it increases the number of faces according to Lemma \ref{lemma:Diff}, which is impossible. Therefore $e_1$ must be a bridge. There is a similar reasoning with $e_2$ which is incident to a single face in $M^{(2)}$. When $e, e_1, e_2$ are bridges, so is $e_3$.

The 4-bond case of Figure \ref{fig:4Bond} cannot occur. Indeed, unhooking an edge $e$ decreases the number of faces by 1. However, after $e$ is unhooked, $e_1$ ($e_2$) has a single incident face in $M^{(1)}$ ($M^{(2)}$), i.e. $|\cI_2(e_1)|\leq 1, |\cI_2(e_2)| \leq 1$ in $M_e$. Unhooking them provides two additional faces according to lemma \ref{lemma:Diff} so that $M$ cannot be dominant.

Finally the case with two 2--bonds displayed in Figure \ref{fig:2BondBis} is not possible either. Say that $e$ and $e_3$ form a 2--bond. Unhooking $e$ leads to $M_e$ which has one face less than $M$. The edges $e$ and $e_3$ are then both bridges in $M_e$, meaning that they are both incident to a single face in both $M^{(1)}$ and $M^{(2)}$. Therefore $e_1$ (or $e_2$) is also incident to a single face in both $M^{(1)}$ and $M^{(2)}$, i.e. $\cI_2(e_1) = \emptyset$ in $M_e$. Unhooking it gives three additional faces, thus leading to more faces than $M$, which is impossible. \qed

\begin{proposition} \label{prop:Black}
If $M\in\cM_{\max}$, two edges incident to a square--vertex and forming a 2--bond are incident on the same black vertex.
\end{proposition}

{\bf Proof.} An edge might either be a bridge, or be part of a 2--bond connecting a black vertex to a square--vertex (i.e. minimal set of two parallel edges), or neither of these two kinds. We denote $p(v)$ the number of the last type of edges incident to the black vertex $v$. They can be part of 2--bonds incident to different square--vertices and $k$--bonds for $k>2$. We now prove by induction that for a dominant $M$, $p(v)=0$ for any black vertex $v$. This way, the edges in $M\in\cM_{\max}$ around black vertices are either bridges or form 2--bonds incident on the same square--vertices.
 
For some black vertex $v$, the $p(v)$ considered edges $(e_1,...,e_{p(v)})$ form an edge-cut since the other edges attached to $v$ are all either bridges or pairs of edges that form 2--bonds. Let $p_a\le p_b\le ...$ be the numbers of edges in each bond of the unique decomposition of the edge-cut $(e_1,...,e_{p(v)})$.

Notice that there are no bridges in the bond decomposition of $(e_1,...,e_{p(v)})$, as we excluded them by definition. Therefore $ p(v)=0$ or $\forall a,\ p_a>1$. In particular, $e_1$ is always part of a $k$--bond with $k>1$ and $p(v)=1$ is impossible.
 
\begin{figure} \centering
\includegraphics[scale=0.35]{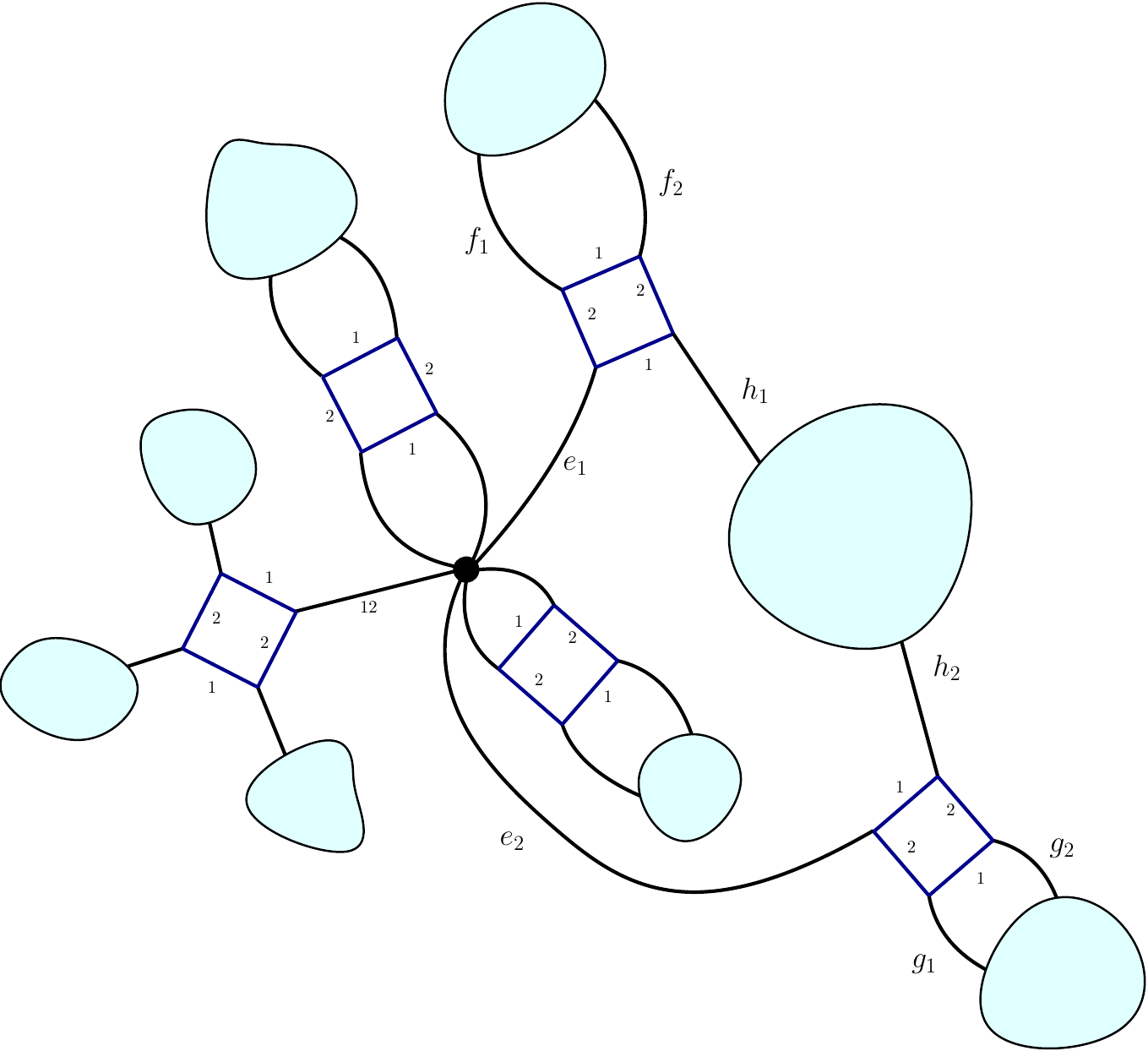}\hspace{0.1cm}\includegraphics[scale=0.37]{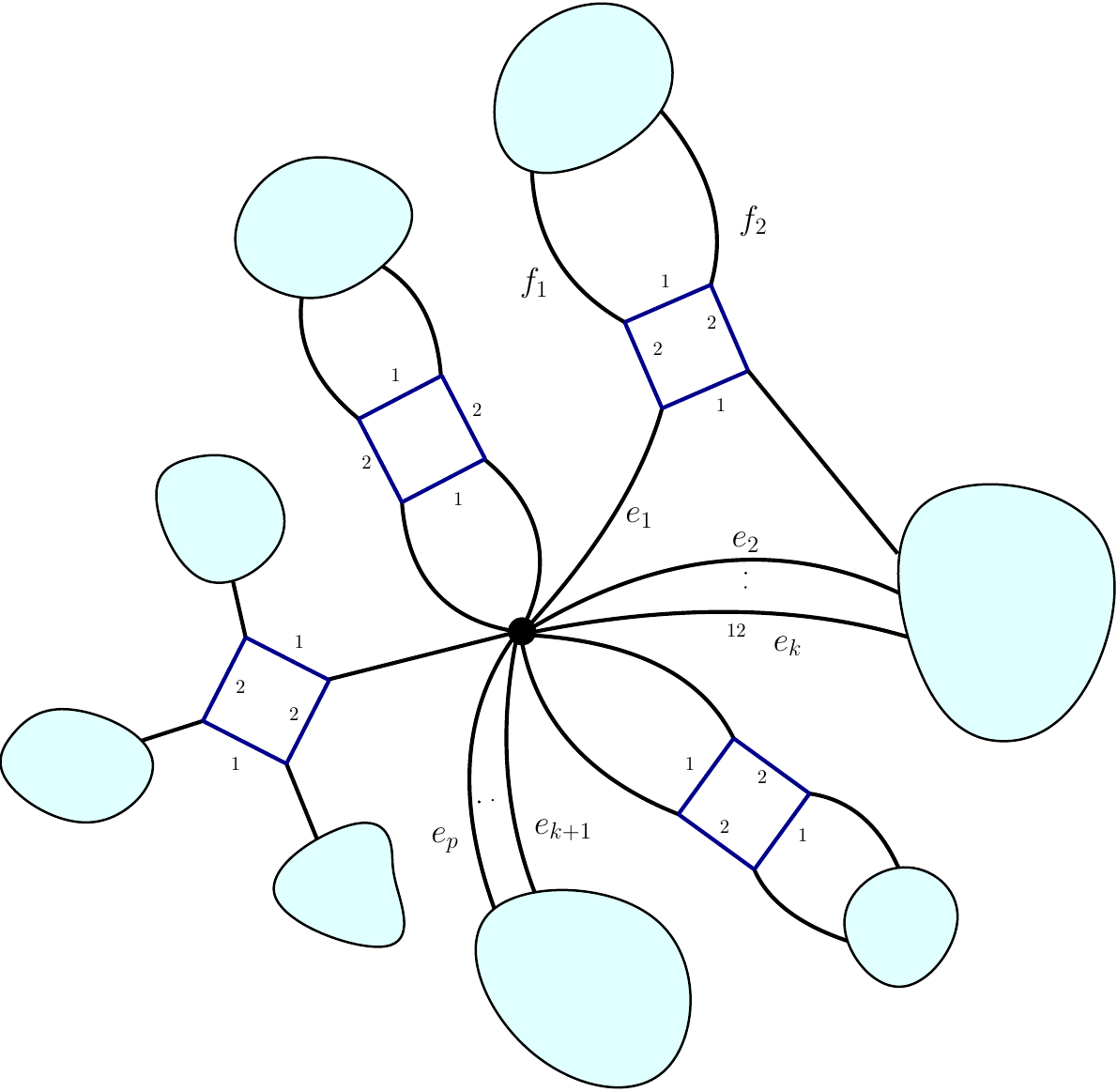}
\caption{\label{fig:Recursion} On the left is a map with $p(v)=2$ and on the right one with $p(v)=k>2$ for some black vertex.}
\end{figure}

Let us look at the case $p(v)=2$, i.e. two edges $e_1$ and $e_2$ which form a 2--bond but reach two different square--vertices. From Proposition \ref{prop:2Bond}, we know that to each of them is attached another 2--bond which does not contain $e_1$ or $e_2$. We denote them $(f_1,f_2)$ and $(g_1,g_2)$ as in the left of Figure of \ref{fig:Recursion} (note that the case in figure \ref{fig:Recursion} is the most general case for which $(e_1,h_1)$ and $(e_2,h_2)$ form two 2--bonds). After unhooking $e_1$ one obtains $M'$ with one face less than $M$ (from lemma \ref{lemma:Diff}) and in which both $e_1$ and $e_2$ are bridges. This implies that $f_1$ and $g_1$ are bridges in at least $M'^{(1)}$ or $M'^{(2)}$, so that unhooking them both brings two additional faces (according to lemma \ref{lemma:Diff}). One thus obtains a new element $M''$ with more faces than $M$, which is impossible.
 
Now suppose that for $q>1$, it is proven that $M\in\cM_{\max}$ has no black vertex $v$ for which $1<p(v)\leq q $, and let $M\in\cM_{\max}$ with $p(v)=q+1$ for some black vertex $v$. Since $e_1$ is not a bridge, $M'$ obtained by unhooking $e_1$ has one face less than $M$. From Proposition \ref{prop:2Bond}, $e_1$ is incident to a square--vertex to which another pair of edges $(f_1, f_2)$ is attached and form a 2--bond as shown on the right of Figure \ref{fig:Recursion}. After detaching $e_1$, we may also detach one of these edges, e.g. $f_1$. As $e_1$ is unhooked, $f_1$ is now a bridge in either $M'^{(1)}$ or $M'^{(2)}$. Let us choose the case $M'^{(2)}$ as in figure \ref{fig:Recursion}. Since $M\in\cM_{\max}$, $f_1$ was not a bridge in $M^{(1)}$ and had two distinct incident faces of color 1, and this is still the case after unhooking $e_1$ (because $e_1$ and $f_1$ belonged to two different connected components in $M^{(1)}$). Unhooking $f_1$ therefore gives a graph $M''$ with one more face than $M'$, hence $F(M'') = F(M)$ and $M''\in\cM_{\max}$. However $M''$ has a vertex $v$ with $p(v)=q$ which contradicts our hypothesis. It is clear that $p(v)=q$ in $M''$ when $q\geq 3$, but the case $q=2$ is more subtle. Indeed, the quantity $p(v)$ counts 2--bonds if their edges are not incident to the same square--vertex and we thus have to rule out that after detaching $e_1$, the two remaining edges of the edge-cut do not form a 2--bond incident on the same square--vertex. This situation corresponds to the case $p(v) = 3$ in $M$ with the edges $e_1, e_2, e_3$ forming a 3--bond. The edges $e_2$ and $e_3$ cannot in fact be incident to the same square--vertex, else they would form a 2--bond. This ensures that $p(v) = 2$ in $M''$. This concludes the induction. \qed

\begin{definition}
Assume that the edges incident to a square--vertex form two 2--bonds like in Figure \ref{fig:2Bond}. The vertical cut of the square--vertex consists in removing the inner edges of color $a$ of the square--vertex which connect both 2--bonds.

Let $M\in\cM_{\max}$ without bridges, and thus satisfying the properties of Propositions \ref{prop:2Bond} and \ref{prop:Black} without bridges. The vertical cut of $M$ is obtained by performing the vertical cut of each square--vertex. It leads to a map $VM$ which has a single black vertex per connected component. We say that $M$ is planar if $VM$ is a planar map.
\end{definition}

\begin{proposition} \label{prop:Planar}
If $M\in\cM_{\max}$ without bridges, then $M$ is planar.
\end{proposition}

{\bf Proof.} From Propositions \ref{prop:2Bond} and \ref{prop:Black}, we know that a square--vertex in $M\in\cM_{\max}$ without bridges is adjacent to exactly two distinct black vertices, and the vertical cut of the square--vertex separates $M$ into two connected components,
\begin{equation} \label{SquareVertexDominant}
\begin{array}{c} \includegraphics[scale=.65]{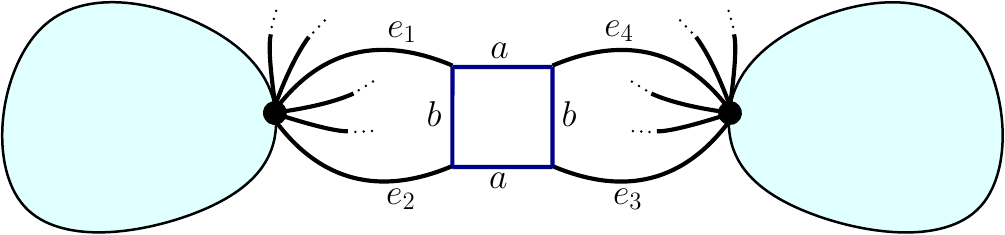} \end{array}
\end{equation}
In the map $M^{(a)}$, one deletes the inner edges of color $b$, and the other way around for $M^{(b)}$. In $M^{(a)}$, the edges $e_1$ and $e_4$ are merged into a single edge, as well as $e_2$ with $e_3$, while in $M^{(b)}$, the edges $e_1$ and $e_2$ are merged into a single loop, as well as $e_3$ with $e_4$. This turns the edges incident to the square--vertex into a pair of parallel edges in $M^{(a)}$ and a pair of loops in $M^{(b)}$,
\begin{equation}
M^{(a)} = \begin{array}{c} \includegraphics[scale=.35]{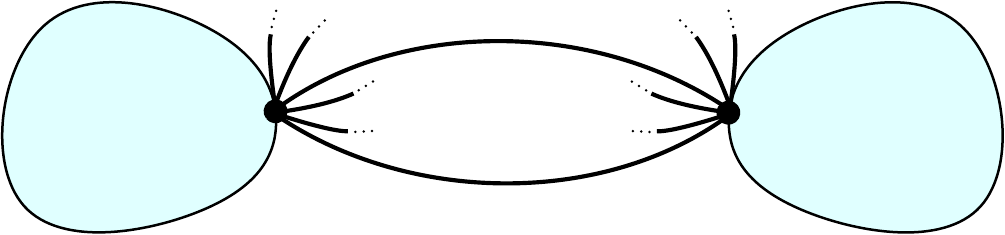} \end{array}, \quad
M^{(b)} = \begin{array}{c} \includegraphics[scale=.35]{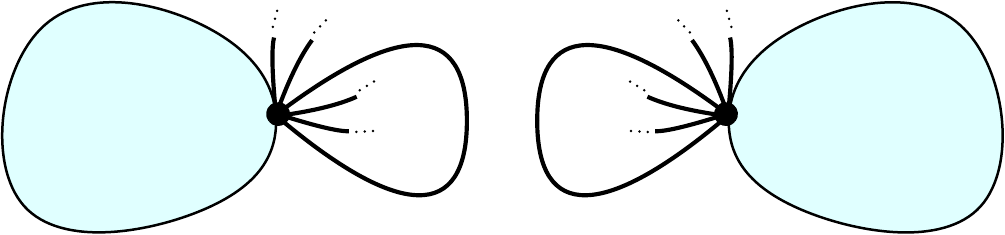} \end{array}
\end{equation}
From Euler's formula, the number of faces of $M^{(c)}$, $c=1, 2$ is
\begin{equation}
F(M^{(c)}) = 2b - V + 2(k^{(c)} - g^{(c)}),
\end{equation}
where we have used the fact that the number of edges of $M^{(c)}$ is $2b$ ($b$ the number of square--vertices) and the number of vertices of $M^{(c)}$ is $V$ the number of black vertices of $M$. Moreover, $k^{(c)}, g^{(c)}$ respectively denote the number of connected components and the genus of $M^{(c)}$.

One can easily turn $M^{(1)}$ and $M^{(2)}$ into planar maps by permuting the order of the edges around the black vertices. For instance, one can make parallel edges occupy consecutive corners. This way, when edges are parallel in $M^{(a)}$, they become a pair of disjoint loops in $M^{(b)}$ (and the other way around) and there is no edges sitting at the corner inside each loop. This permuting of edges around black vertices does not change the number of connected components of $M^{(c)}$, but only its genus. This thus maximizes $F(M^{(c)})$ and it can be concluded that $M\in\cM_{\max}$ without bridges has planar maps $M^{(1)}$ and $M^{(2)}$.

Let $v$ be a black vertex in $M$. The faces of $M^{(c)}$, $c=1, 2$, can be partitioned as those which go through $v$ (there are $F^{(c)}_v$ of them) and those which do not (there are $F^{(c)}_{\hat{v}}$ of them),
\begin{equation}
F(M^{(c)}) = F^{(c)}_v + F^{(c)}_{\hat{v}}
\end{equation}
Since $M^{(1)}$ and $M^{(2)}$ are planar, there are well defined notions of outside and inside the faces which are delimited by either parallel edges or loops. The face which leaves the black vertex $v$ along the external (resp. internal) side of an edge which is part of a 2--bond in $M^{(a)}$ comes back to $v$ for the first time along the external (resp. internal) side of the edge in the same 2--bond,
\begin{equation}
\begin{array}{c} \includegraphics[scale=.5]{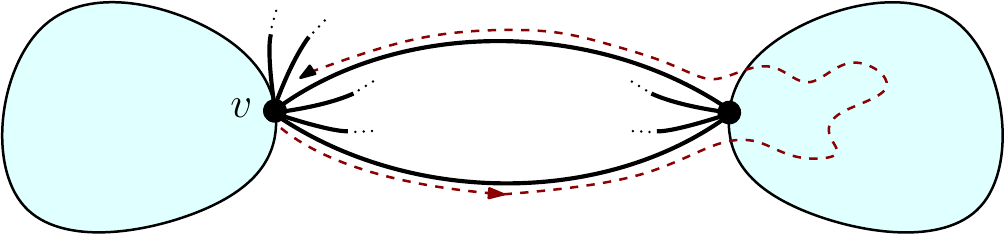} \end{array}
\end{equation}
This property is obvious for $M^{(b)}$ as one travels outside or inside a loop. With some abuse of notation, it can be visualized on $M$ itself, by saying that the external (internal) face leaving $v$ along $e_2$ returns to $v$ for the first time along the external (internal) side of $e_1$,
\begin{equation}
\begin{array}{c} \includegraphics[scale=.5]{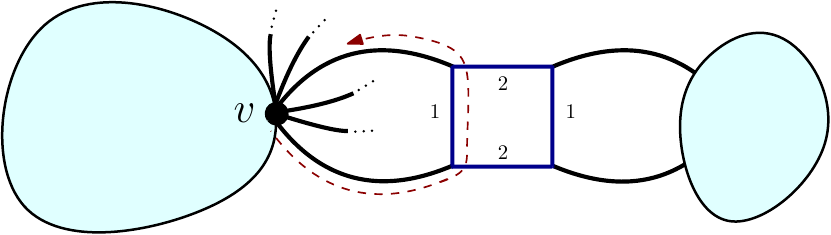} \end{array} \qquad \begin{array}{c} \includegraphics[scale=.5]{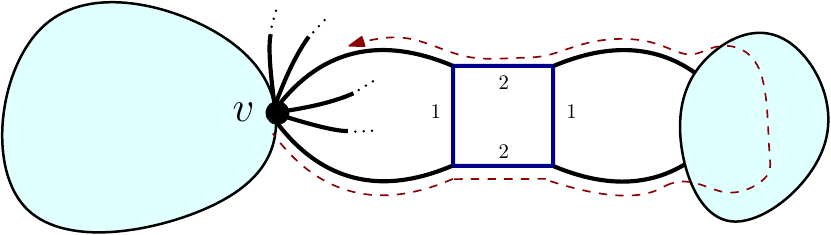} \end{array}
\end{equation}
This shows that there is a bijection between the faces which go through $v$ in $M^{(1)}$ and $M^{(2)}$, and $F^{(1)}_v = F^{(2)}_v$.

Let $M_0(v)$ be the connected component of the vertical cut of $M$ which contains $v$. For the square--vertex represented in \eqref{SquareVertexDominant}, one deletes the inner edges of color $a$ and gets two loops. In fact each square--vertex provides the vertical cut with two loops. Therefore, there is a bijection between the faces of $M_0(v)$ and the faces of $M^{(1)}$ which go through $v$ (and $M^{(2)}$ as well) and thus
\begin{equation}
F(M^{(c)}) = F(M_0(v)) + F^{(c)}_{\hat{v}}
\end{equation}
At fixed $F^{(c)}_{\hat{v}}$, one maximizes $F(M^{(c)})$ by maximizing $F(M_0(v))$. Since $M_0(v)$ is a 1--vertex map, this is done by selecting any planar configuration for $M_0(v)$. This reasoning applies to any black vertex of $M$ and shows that the vertical cut has to be planar. \qed

\begin{theorem} \label{thm:Dominant}
The set $\cM_{\max}$ is defined by the Propositions \ref{prop:2Bond}, \ref{prop:Black} and \ref{prop:Planar}. An element $M\in\cM_{\max}$ with $b(M)$ square--vertices has $F(M)$ faces given by $F(M) = 5b(M) + 3$.
\end{theorem}

\begin{figure} \centering
\includegraphics[scale=.3]{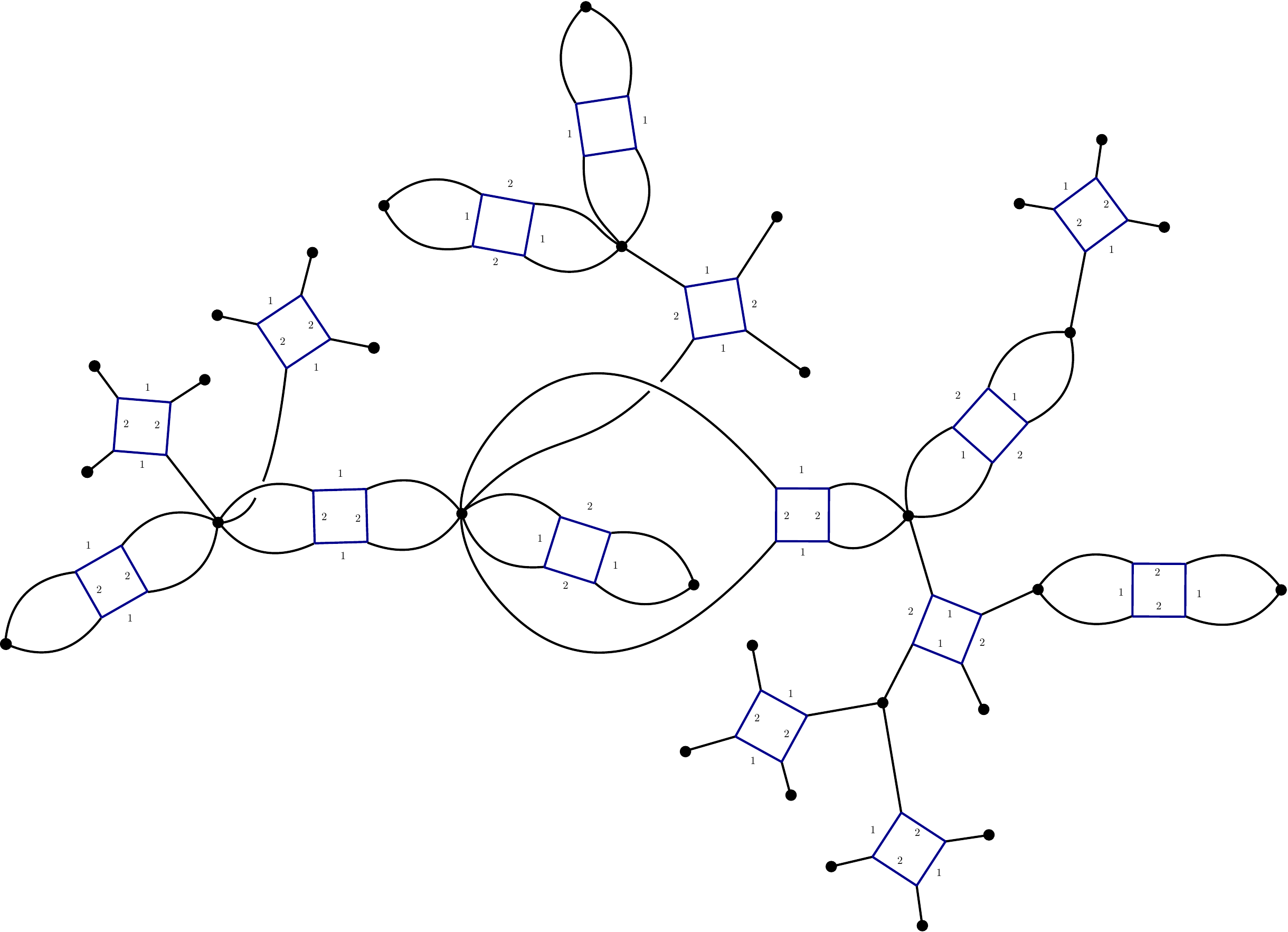}
\caption{\label{fig:Dominant} A generic example of a dominant map.}
\end{figure}

An example of a generic dominant map is given in figure \ref{fig:Dominant}.

{\bf Proof.} Let $\cM_0$ be the set of elements of $\cM$ satisfying the criteria given in Propositions \ref{prop:2Bond}, \ref{prop:Black}, \ref{prop:Planar}. From those propositions, we already know $\cM_{\max} \subset \cM_0$. It is therefore sufficient to show that all elements of $\cM_0$ have $5b(M) + 3$ faces.

The criteria of Propositions \ref{prop:2Bond}, \ref{prop:Black}, \ref{prop:Planar} show that around a square--vertex of $M\in\cM_0$, $M$ takes any one of the following forms
\begin{equation}
\begin{array}{c} \includegraphics[scale=.45]{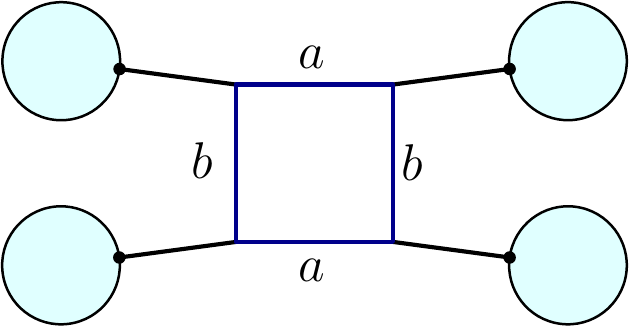} \end{array} \qquad \begin{array}{c} \includegraphics[scale=.45]{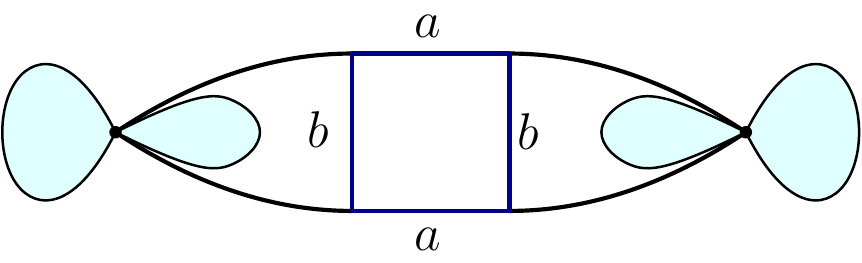} \end{array}
\end{equation}
where the blobs reproduce the same pattern. In the case where the square--vertex is incident to two 2--bonds, one can unhook an edge in each 2--bond. Unhooking one decreases the number of faces by one. After that, the edges of the other 2--bond each have $|\cI_2(e)|=1$ (they are bridges in the color $a$ but still have two incident faces in the color $b$). Therefore, one gets from lemma \ref{lemma:Diff} that unhooking a second edge gives a new element of $\cM$ with exactly the same number of faces as $M$. Doing so for all square--vertices incident to 2--bonds, the number of faces is preserved and $M$ is transformed into a tree for which it is known from Proposition \ref{prop:Trees} that $F(M) = 5b(M) + 3$. \qed

\subsection{The colored graphs behind the dominant gluings} \label{sec:DominantColoredGraphs}

It is clear from the above analysis that if $\cM$ is in $\cM_{\max}$ with $b$ square--vertices, then it has a square--vertex which is incident to only one other and looks like one of the two possibilities below,
\begin{equation}
\begin{array}{c} \includegraphics[scale=.6]{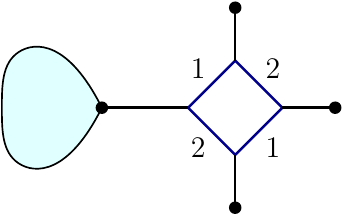} \end{array} \hspace{4cm} \begin{array}{c} \includegraphics[scale=.6]{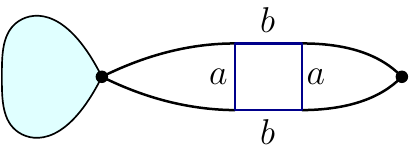} \end{array}
\end{equation}
where the blobs are elements of $\cM_{\max}$ with $b-1$ square--vertices. From the point of view of the blobs, the case on the left is the insertion of a square--vertex via a bridge, and on the right via a 2--bond, with $a = 1, 2$.

Those insertions can be performed on arbitrary corners around black vertices and always create elements of $\cM_{\max}$. From the point of view of colored graphs, a corner is an edge of color 0. The above bridge or 2--bond insertions correspond to cutting the edge of color 0 into two halves and gluing those halves onto a new bubble. Since the square--vertices with four incident bridges and with two incident 2--bonds are equivalent to the colored graphs given in \eqref{CubePairing}, we find that the elements of $\cM_{\max}$ are in bijection with the elements of $\cG_{\max}\subset \cG$ which we now describe.

Denote $P_1, P_2, P_3$ the graphs in \eqref{CubePairing} where the edges of color 0 are parallel to those of color $a$ in $P_a$. An element $G_b$ of $\cG_{\max}$ with $b$ bubbles is always obtained from an element $G_{b-1}$ of $\cG_{\max}$ with $b-1$ bubbles by cutting open an edge of color 0 of $G_{b-1}$, an edge of color 0 of $P_a$ for $a\in\{1, 2, 3\}$, and reconnecting the half-edges to get a connected, bipartite graph,
\begin{equation}
\begin{array}{c} \includegraphics[scale=.6]{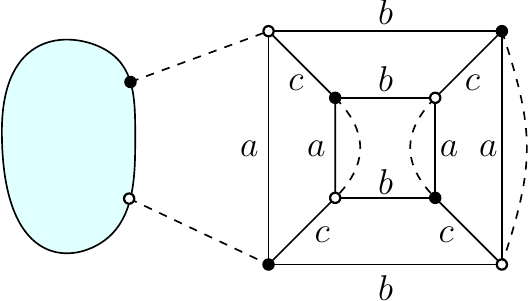} \end{array}
\end{equation}

For instance, the graphs with two bubbles are obtained by taking $P_a$ and $P_b$ and cutting an edge parallel to the color $a$ in $P_a$, an edge parallel to the color $b$ in $P_b$ and reconnecting them to connect $P_a$ with $P_b$.


\section{Topology} \label{sec:Topology}

Although we are not primarily interested in the topology of the gluings, we can easily show that elements of $\cM_{\max}$ correspond to colored triangulations of the 3--sphere.

First, one shows that the three triangulations, represented in \eqref{CubePairing}, with a single octahedron and which maximize the number of edges, are spheres. From the symmetry of the colors, it is sufficient to prove it for only one of them, say the one on the left of \eqref{CubePairing}. To do so, notice that if one provides the colored graph with a well--chosen cyclic order of the colors around its vertices, then it becomes a planar map. Indeed, choose the cyclic order $(0 1 3 2)$ around white vertices and its reverse cycle $(0 2 3 1)$ around black vertices, as actually depicted in \eqref{CubePairing}. This cyclic order makes it a map, called a jacket, which turns out to be planar. It is known in the topological theory of colored graphs that if a graph with four colors (i.e. in three dimensions) admits a planar jacket, then it is a sphere (see for instance \cite{GurauRyanReview}).

This paves the way for an induction. Consider $M\in\cM_{\max}$ and via the bijection of Section \ref{sec:Bijection} let $G$ be the corresponding graph. If $M$ has at least two square--vertices (i.e. $G$ has at least two bubbles), then it has a black vertex which is a cut--vertex (from the planarity of Proposition \ref{prop:Planar}). In terms of colored graphs, this black vertex corresponds to a cycle of $G$ alternating edges of colors 0 and 3. From the planarity of Proposition \ref{prop:Planar}, this cycle has a pair of edges of color 0 forming an edge--cut,
\begin{equation}
\begin{array}{c} \includegraphics[scale=.4]{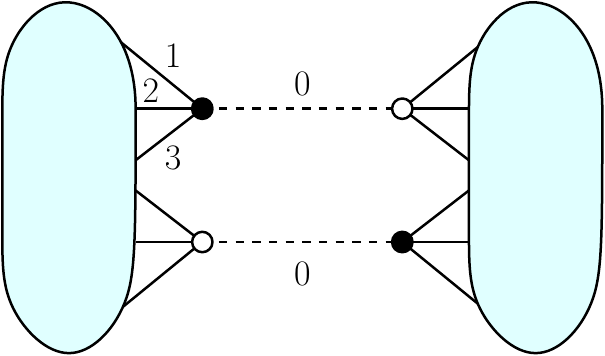} \end{array}
\end{equation}
Performing the cut and reconnecting the edges of color 0 on the left and on the right, one obtains two connected components,
\begin{equation}
\begin{array}{c} \includegraphics[scale=.4]{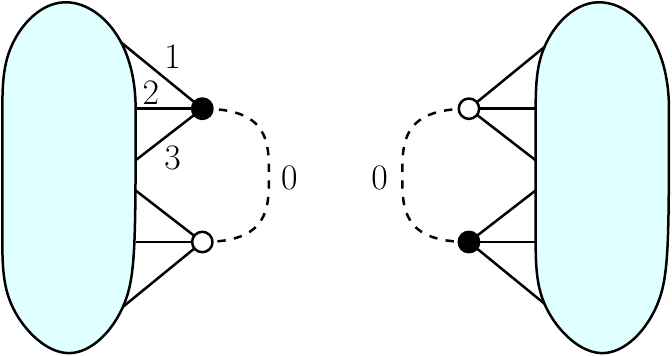} \end{array}
\end{equation}
They are both graphs with less bubbles than $G$ and which maximize the number of bicolored cycles for their numbers of bubbles. Let us assume that they are both spheres. We can now reconstruct $G$ by the following sequence of operations. First, one inserts 3--dipoles on the two edges of color 0,
\begin{equation}
\begin{array}{c} \includegraphics[scale=.4]{2EdgeCutDisconnected.pdf} \end{array} \qquad \to \qquad \begin{array}{c} \includegraphics[scale=.4]{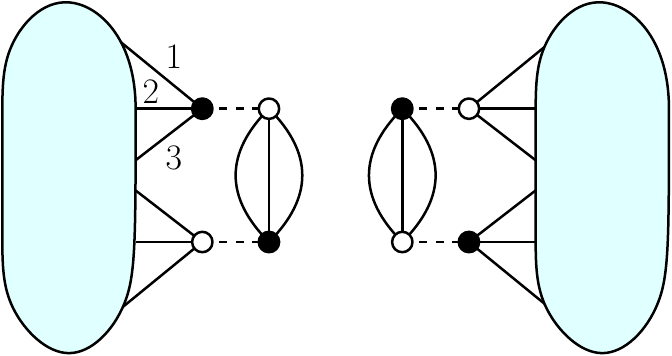} \end{array}
\end{equation}
A 3--dipole is a pair of vertices connected by all the colors except 0. Topologically, it is a 3--ball consisting of two tetrahedra glued along three triangles and whose boundary consists of two triangles of color 0. The insertion of a 3--dipole therefore does not change the topology. Then, we contract a black vertex with a white vertex on both sides. It means that the two chosen vertices are removed and their incident edges are attached pairwise respecting the colors. We choose to do it on vertices of the 3--dipoles (here circled by a dotted line),
\begin{equation}
\begin{array}{c} \includegraphics[scale=.4]{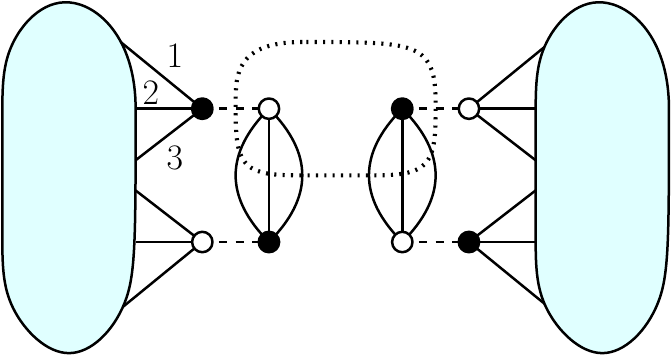} \end{array} \qquad \to \qquad \begin{array}{c} \includegraphics[scale=.4]{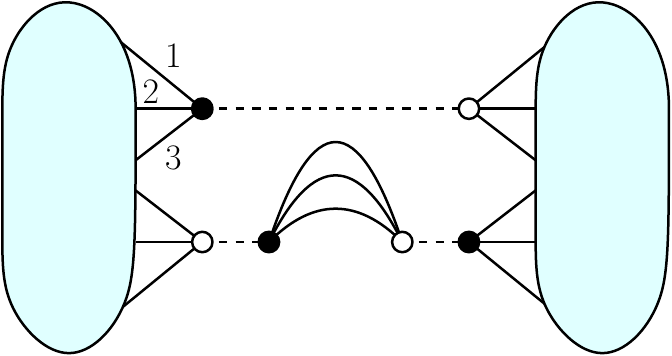} \end{array}
\end{equation}
This contraction operation is topologically a connected sum \cite{ItalianSurvey}, so that it leads to a sphere. Finally, the graph obtained is simply $G$ with an additional 3--dipole on one edge of color 0. Removing that 3--dipole gives $G$ without topology changes. Therefore $G$ is a sphere too.

\section{Enumeration} \label{sec:Enumeration}

We now consider elements of $\cM_{\max}$ equipped with a root, i.e. a distinguished corner at a black vertex called the root vertex. Due to the cyclic order at each black vertex, marking a corner is equivalent to distinguishing an edge (the one following the marked corner clockwise, for instance), which we call the root edge. Let $\cA$ be the set of rooted elements of $\cM_{\max}$ and $A(z) = \sum_{M\in \cA} z^{b(M)}$ its generating function with respect to the number of square--vertices.

\subsection{Using edge deletion and root vertex degree} \label{sec:EnumerationEdgeDeletion}

It is standard in map enumeration to consider the process of deleting the root edge in order to obtain a recursion. We will follow this method in this subsection. This usually necessitates to take into account the degree of the root vertex. Let $\cA_n$ be the set of rooted elements of $\cM_{\max}$ where the root vertex has degree $n$. $\cA_0$ has one element which consists of an isolated black vertex. The set of rooted elements of $\cM_{\max}$ decomposes as $\cA = \sqcup_{n\geq 0} \cA_n$, and we have the following generating functions 
\begin{equation}
A_n(z) = \sum_{M\in\cA_n} z^{b(M)}, \qquad A(z) = \sum_{M\in \cA} z^{b(M)} = \sum_{n\geq 0} A_n(z).
\end{equation}

The root edge connects a black vertex to a square--vertex. According to Proposition \ref{prop:2Bond}, it is either a bridge or part of a 2--bond, which induces the partition
\begin{equation}
\cA_n = \cB_n \sqcup \cC_n
\end{equation}
where $\cB_n$ (resp. $\cC_n$) are rooted elements of $\cM_{\max}$ with a root edge on a black vertex of degree $n$ which is a bridge (resp. part of a 2--bond). We denote their generating functions $B_n(z)$ and $C_n(z)$ respectively.

\begin{enumerate}
\item \label{enum:Bridge} If the root edge $e$ is a bridge, then it is incident to a black vertex $v$ which has $n-1$ other incident edges. The latter can form any contribution $M_0$ from $\cM_{\max}$ and rooted, e.g. on the corner incident to $e$ at $v$ in the clockwise direction, hence $M_0\in\cA_{n-1}$. We further know from Proposition \ref{prop:2Bond} that the three other edges incident to the same square--vertex as $e$, say $e_1, e_2, e_3$, are bridges too. They are moreover connected to distinct black vertices to which any element of $\cA$ can be attached. Indeed, those contributions are obviously arbitrary elements of $\cM_{\max}$ but they are also rooted, e.g. on the corners sitting clockwise to $e_1, e_2, e_3$, and there are no restrictions on the degrees of the vertices to which $e_1, e_2, e_3$ are attached. This leads to
\begin{equation}
B_n(z) = z\,A_{n-1}(z)\,A(z)^3. 
\end{equation}
In particular, the generating function of $\cA_1$ satisifies
\begin{equation}
A_1(z) = z\,A(z)^3
\end{equation}
since $\cA_1 = \cB_1$ and $A_0(z)=1$.

\item \label{enum:NonBridge} If the root edge $e$ is not a bridge, it is part of a 2--bond. Denote $(e, e_1, \dotsc, e_{n-1})$ the edges incident to $v$, cyclically ordered. The root edge $e$ forms a 2--bond with $e_{k+1}$ for $k\in\{0, \dotsc, n-2\}$. They are connected to a square--vertex, in turn incident to a 2--bond $\{f, g\}$ incident on a black vertex $v'$. The edges $e, e_{k+1}$ divide the neighborhood of $v$ into two angular regions, and similarly with $f, g$ at $v'$. Planarity from Proposition \ref{prop:Planar} imposes that all contributions to $M$ are attached to $v$ or $v'$ within those angular regions. This gives four elements of $\cA$, say $M_1$ rooted on the corner next to $e$ in the clockwise direction, $M_2$ rooted on the corner next to $e_k$ in the clockwise direction, $M_3$ rooted on the corner next to $f$ in the clockwise direction and $M_4$ rooted on the corner next to $g$ in the clockwise direction. Moreover $M_1$ is rooted on a vertex with $k$ edges $(e_1, \dotsc, e_k)$ and $M_2$ on a vertex with $n-k-2$ edges $(e_{k+2}, \dotsc, e_{n-1})$, i.e. $M_1\in\cA_{k}$ and $M_2\in\cA_{n-k-2}$ while there are no degree restrictions on the root of $M_3$ and $M_4$. Finally there are two ways to position the square--vertex so that $\{e, e_{k+1}\}$ is a 2--bond (with an inner edge of color 1 or of color 2 between them). This way it is found that
\begin{equation}
C_n(z) = 2 z\,A(z)^2 \sum_{k=0}^{n-2} A_k(z)\,A_{n-k-2}(z).
\end{equation}
\end{enumerate}

From $A_n(z) = B_n(z) + C_n(z)$ for $n\geq 2$, one gets
\begin{equation} \label{RecursionA_n}
A_n(z) = z\,A(z)^3\,A_{n-1}(z) + 2 z\,A(z)^2 \sum_{k=0}^{n-2} A_k(z)\,A_{n-k-2}(z),
\end{equation}
for $n\geq 2$ and $A_1(z) = z A(z)^3$ and $A_0(z) = 1$. We sum those equations for all $n\geq 2$ so as to get an equation on $A(z) = \sum_{n\geq 0} A_n(z)$, which remarkably simplifies to
\begin{equation} \label{CountingEq}
A(z) = 1 + 3z\,A(z)^4.
\end{equation}
Expanding $A(z)$ as a power series $A(z) = \sum_{k\geq 0} a_k z^k$, this is equivalent to the recursion
\begin{equation}
a_{n+1} = 3 \sum_{\substack{k_1, k_2, k_3, k_4\\ k_1+k_2+k_3+k_4 = n}} a_{k_1}a_{k_2}a_{k_3}a_{k_4}
\end{equation}
with $a_0=1$. One finds $(1, 3 , 36, 594, \dotsc)$.

The singularity analysis of $A(z)$ is completely straightforward. Singular points $(A_c, z_c)$ are solutions of $\Phi(A, z) \equiv 1 - A + 3z A^4 = 0$ and $\partial_A \Phi(A, z) = -1 + 12z A^3 = 0$. The second equation gives $z_c = 1/(12 A_c^3)$ which after being plugged into the first equation leads to $A_c = 4/3$ and $z_c = 9/256$. Expanding $1-A(z)+3zA(z)^4 = 0$ around that point gives the critical behavior
\begin{equation}
A(z) = \frac{4}{3} - \sqrt{\frac{2048}{243} \Bigl(\frac{9}{256} - z\Bigr)} + o\Bigl(\sqrt{\frac{9}{256} - z}\Bigr).
\end{equation}
which lies in the universality class of trees.

\subsection{Direct enumeration} \label{sec:DirectCounting}

Equation \eqref{CountingEq} suggests that the enumeration should be similar to a tree counting and does not require to take into account the degree of the root vertex. This can also be seen from the fact that no problem was encountered when summing the equations \eqref{RecursionA_n} over $n$ without using a catalytic variable, in contrast to the generic case of planar maps where this is not directly possible.

A direct counting\footnote{We would like to thank an anonymous referee for pointing out this direct counting to us.} can be given which clearly shows the tree structure of the maps in $\cA$. Starting from the root corner and following the cyclic order around the root vertex (say clockwise), one can encounter an arbitrary sequence of two types of contributions. They correspond to the two cases of the above counting: they repeat the same arguments, simply without taking into account the degree of the root vertex (and taking a sequence of them). We re-state them quickly:
\begin{itemize}
\item Case \ref{enum:Bridge}. One can encounter a bridge connected to a square--vertex, itself connected to three other bridges. Those three other bridges are attached to arbitrary rooted maps in $\cA$, thus contributing like $z A(z)^3$.

\item Case \ref{enum:NonBridge}. One can encounter an edge forming a 2--bond with another edge connected to the root vertex. Between those two edges, any rooted element of $\cA$ can be found, contributing to $A(z)$. They are moreover connected to a square--vertex to which another pair of edges is attached. The latter form a 2--bond and are connected to another black vertex. They delimit two regions around this black vertex on which two arbitrary rooted maps can be found, contributing to $A(z)^2$. This case thus contributes like $2z A(z)^3$, where the factor $2$ is due to the two choices of labels on the square--vertex.
\end{itemize}
Taking a sequence of those two types of contributions leads to 
\begin{equation} \label{PlaneTrees}
A(z) = \sum_{n\geq 0} \bigl(3z A(z)^3)^n = \frac{1}{1-3z  A(z)^3},
\end{equation}
as desired.

\subsection{Explicit bijection with plane trees}

Equation \eqref{CountingEq} suggests a bijection with trees. A direct interpretation would be as quarternary trees, but one can also write the equation as $A(z) = 1/(1-3z A(z)^3)$ which suggests an interpretation as plane trees of arbitrary degrees. The latter case is more directly found. 

One way to do it is inspired by the direct counting provided above in Section \ref{sec:DirectCounting}. That counting is based on a specification of elements of $\cA$ which clearly unravels the tree--like structure. In this section, we make it more explicit by constructing a bijection with plane trees. Since we feel that it is quite straightforward from the previous section of the article and it involves classical combinatorial steps, we only sketch the construction.

First, we notice that the root induces a root (i.e. a marked corner) on each black vertex. Indeed, start from the root edge. If it is a bridge, we have already explained above how to root the three other black vertices adjacent to the same square--vertex (Case \ref{enum:Bridge} of Section \ref{sec:EnumerationEdgeDeletion}). If it is part of a 2--bond, there is another 2--bond attached to the same square--vertex and we have already explained above how to root the black vertex of that 2--bond (Case \ref{enum:NonBridge} of Section \ref{sec:EnumerationEdgeDeletion}). Next, we repeat the same procedure from these new marked corners. Furthermore, if the second edge sitting next to the root edge clockwise is not the second edge of the 2--bond containing the root edge, we also apply this procedure as if the second edge were the root edge, and so on. This way, every black vertex is visited exactly once and receives a marked corner induced by the root.

The next step is to associate a tree to the neighborhood of each black vertex. There are several rather standard ways to do so. Here are two constructions which lead to the same tree.
\begin{itemize}
\item Consider the vertical cut of $M\in\cA$, obtained by deleting the inner edges of square--vertices which separate 2--bonds, and delete the square--vertices incident to bridges. Assume for a moment that there are no bridges. Then the vertical cut has one connected component per black vertex, by construction, and each of them is thus a rooted planar 1--vertex map. It is well known that they are in bijection with planar trees: this bijection is obtained by taking the dual of the rooted planar 1--vertex map. Each edge of such a tree corresponds to a 2--bond of $M$. Since there are two types of 2--bonds (depending on the color of the inner edges adjacent to the edges of the 2--bond), the edges of the tree can be colored with two possible colors. Moreover, the bridges can be added back into corners of the 1--vertex maps. Since those corners are mapped to corners of the dual trees, bridges can be inserted into the corners of the trees, and they have to be given a new color to remember that they come from bridges in $M$.
\item An alternative is to visit each edge clockwise around a black vertex. Notice that the edges forming 2--bonds come in canonical pairs and for each of them there is an opening edge and a closing edge determined by the clockwise order starting from the root corner. We say that an edge is the child of the previous one if it is immediately after the opening edge of a 2--bond. We say that it is a sibling of the previous edge if the latter is either a bridge or the closing edge of a 2--bond. These children/sibling relationships give rise to a plane tree for each black vertex. Moreover, an edge of the tree corresponding to one of a 2--bond can have two colors (for the two choices of labels of the square--vertex), and an edge of tree corresponding to a bridge receives a third color.
\end{itemize}
One thus obtains for $M$ a collection of rooted planar trees whose edges can have three possible colors (note that if the color indicates a bridge in $M$, this stops the branch). This is still not enough to fully characterize $M$.

Consider two 2--bonds in $M$ incident to the same square--vertex. They are mapped to two edges in two different trees. To reconstruct the 2--bonds of $M$ from the trees, it is sufficient to pair (with a dashed edge connecting them for instance) the edges of the trees which corresponded to 2--bonds incident to the same square--vertex. In case of bridges, there are four edges in four different trees which have to be connected by dashed edges to remember that they are incident to the same square--vertex of $M$.

We can now repeat the counting of those objects. The root vertex (incident to the root edge) can have an arbitrary (finite) number $n\geq 0$ of children. In terms of generating function, the case $n=0$ gives 1. For $n>0$, the edges between the root and its children have three possible colors. If the color indicates a bridge in $M$, then this stops the branch. The edges are however connected via dashed edges to three other edges with the same color. Each of them has a (rooted planar) parent tree. Therefore an edge incident to the root with the color associated to bridges of $M$ gives a factor $zA(z)^3$ to the generating function.

In the case where the edge has a color indicating that it is a 2--bond in $M$, then it can have an arbitrary (rooted planar) child tree, hence contributing $A(z)$ to the generating function. It is moreover connected via a dashed edge to another edge representing a 2--bond. This edge in another tree has an arbitrary (rooted planar) parent tree and an arbitrary (rooted planar) child tree, yielding a factor $A(z)^2$. Since there are two colors in this case, the generating function receives a factor $2zA(z)^3$. Altogether, the counting equation reads exactly \eqref{PlaneTrees} as expected.

\section{Conclusion}

In this article, we have performed the first analysis of colored triangulations whose building blocks are non-melonic and have spherical topology, using the octahedron. The main challenge is that for non-melonic building blocks, Gurau's degree \eqref{GurauDegree} never vanishes and it is expected that there is no infinite families at constant values of the degree. Therefore we had to identify an enhanced degree found to be
\begin{equation}
\tilde{\omega}_B(T) = \frac{5}{8} t(T) + 3 - e_0(T)
\end{equation}
where the factor $3/4$ of Gurau's degree for bubbles with eight tetrahedra \eqref{GurauDegreeBubble} is replaced with $5/8<3/2$. We further identified the gluings of octahedra which maximize the number of edges, corresponding to $\tilde{\omega}_B(T) = 0$ and showed that they are in bijection with a family of trees and can be enumerated explicitly.

This is similar to the case of arbitrary colored triangulations (with unconstrained building blocks): if one looks for the colored triangulations which maximize the number of edges among the set of all colored triangulations, one finds the melonic triangulations which are in bijection with trees. This is the premise of the classification of colored triangulations according to Gurau's degree performed in \cite{GurauSchaeffer}. This suggests to apply the techniques of \cite{GurauSchaeffer} to classify the gluings of octahedra according to the enhanced degree $\tilde{\omega}_B(T)$.

We showed in Section \ref{sec:Enumeration} that the generating function of dominant graphs with a marked edge of color 0 satisfies $A(z) = 1 + 3z A(z)^4$. In fact, the arguments of Sections \ref{sec:Enumeration} and \ref{sec:DominantColoredGraphs} show that the generating function $P(z)$ of dominant graphs with a marked bubble is
\begin{equation} \label{CubeExpectation}
P(z) = 3z A(z)^4.
\end{equation}
The interpretation is given in Section \ref{sec:DominantColoredGraphs}. One starts with the marked bubble. There are three ways to add four edges of color 0 and maximize the number of bicolored cycles with this bubble only, shown in \eqref{CubePairing}. Now to get $P(z)$, one can cut any of those edges of color 0 and add a dominant graph with a marked edge of color 0 instead, whose generating function is $A(z)$.

It would be interesting to study the generating functions of graphs with a marked bubble chosen arbitrarily. In the case where dominant graphs are defined by the vanishing of Gurau's degree \eqref{GurauDegreeBubble} (instead of the vanishing of \eqref{OctahedraDegree} with $s_B = 5/8$ in our case), it is known \cite{Universality} that the generating function $P_{B'}(z)$ of dominant graphs with an arbitrary marked bubble $B'$ is
\begin{equation} \label{GaussianUniversality}
P_{B'}(z) = C_{B'} z \tilde{A}(z)^{p_{B'}}
\end{equation}
where $p_{B'}$ is the number of black vertices of $B'$, $\tilde{A}(z)$ the generating function of graphs of vanishing Gurau's degree with a marked edge of color $0$ and $C_{B'}$ is the number of ways to add $p_B$ edges of color 0 to $B'$ and maximize the number of bicolored cycles. Obviously, our calculation of $P(z)$, \eqref{CubeExpectation}, is reminiscent of \eqref{GaussianUniversality}. However, the switch from Gurau's degree to \eqref{OctahedraDegree} is {\it a priori} non--trivial and the theorem of \cite{Universality} cannot be applied in our case. The result \eqref{CubeExpectation} still suggests that there might exist an extension of \cite{Universality} to our case.

The result of \cite{Universality} was framed in the context of random tensors. Indeed, colored triangulations can be generated by specific probability distributions over a random tensor of size $N\times N\times N$. A colored triangulation made of gluings of a bubble $B$ and generated by a random tensor model is moreover weighted by $N^{-\omega_B(T)}$ where $\omega_B(T)$ is Gurau's degree (adapted to $B$) \eqref{GurauDegreeBubble}. In this context, the theorem \eqref{GaussianUniversality} of \cite{Universality} states that the distribution on the random tensor becomes Gaussian when $N\to\infty$. In our case, i.e. $B$ corresponding to octahedra and replacing Gurau's degree with \eqref{OctahedraDegree}, there still exists a random tensor model which generates our gluings of octahedra weighted by $N^{-\tilde{\omega}_B(T)}$. The large $N$ limit therefore also restricts to dominant graphs (defined by $\tilde{\omega}_B(T) = 0$). It would thus be very interesting to know if that distribution becomes Gaussian when $N\to\infty$.

The case of octahedral gluings thus falls into the same universality class as the melonic one. However, the analysis is different. In the case of melonic bubbles, the most efficient strategy is to find a bijection in which the triangulations which maximize the number of edges are just (plane) trees (see for instance \cite{DartoisDoubleScaling, MelonoPlanar}). Then, given an arbitrary triangulation in this representation, one can take a spanning tree in it and show that deleting the edges which are not in the spanning tree corresponds to increasing the number of edges of the corresponding triangulation. In the octahedral case, we have used a representation of the colored triangulations which generalize that of \cite{DartoisDoubleScaling, MelonoPlanar} and which is a special case of \cite{StuffedWalshMaps}. However, the triangulations which maximize the number of edges are not simply trees in this representation, which makes the analysis a bit less straightforward. It would of course be very interesting to find a bijection for the gluings of octahedra such that those maximizing the number of edges directly become trees. If such a bijection could be found, it would be interesting to try and extend it to more generic bubbles in three dimensions. This way, it would actually be a promising strategy to investigate whether other universality classes can be found.

\section*{Acknowledgements}

This research was supported by the ANR MetACOnc project ANR-15-CE40-0014.



\begin{thebibliography}{99}

\bibitem{CatalyticVariables}
  M.~Bousquet-M\'elou and A.~Jehanne,
  ``Polynomial equations with one catalytic variable, algebraic series and map enumeration,''
  J. Combin. Theory Ser. B, 96(5):623-672, 2006.
  \doi{10.1016/j.jctb.2005.12.003}

\bibitem{MatrixReview}
  P.~Di Francesco, P.~H.~Ginsparg and J.~Zinn-Justin,
  ``2-D Gravity and random matrices,''
  Phys.\ Rept.\  {\bf 254}, 1 (1995)
  \doi{10.1016/0370-1573(94)00084-G},
  \arxiv{hep-th/9306153}.

\bibitem{Schaeffer}
  G.~Schaeffer,
  ``Bijective census and random generation of Eulerian planar maps with prescribed vertex degrees,''
  Electron. J. Combin. 4 (1997) R20.

\bibitem{Mobiles}
  J.~Bouttier, P.~Di~Francesco and E.~Guitter,
  ``Planar maps as labeled mobiles,''
  Electron. J. Combin. 11 (2004) R69
  \arxiv{math/0405099}.

\bibitem{TR}
  B.~Eynard,
  ``Topological expansion for the 1-Hermitian matrix model correlation functions,''
  JHEP {\bf 0411}, 031 (2004)
  \doi{10.1088/1126-6708/2004/11/031},
  \arxiv{hep-th/0407261}.

\bibitem{GurauRyanReview}
  R.~Gurau and J.~P.~Ryan,
  ``Colored Tensor Models - a review,''
  SIGMA {\bf 8}, 020 (2012)
  \doi{10.3842/SIGMA.2012.020},
  \arxiv{1109.4812} [hep-th].

\bibitem{SigmaReview}
  V.~Bonzom,
  ``Large $N$ limits in tensor models: Towards more universality classes of colored triangulations in dimension $d\geq 2$,''
  SIGMA {\bf 12}, 073 (2016)
  \doi{10.3842/SIGMA.2016.073},
  \arxiv{1603.03570} [math-ph].

\bibitem{Tensors}
  J.~Ambjorn, B.~Durhuus and T.~Jonsson,
  ``Three-Dimensional Simplicial Quantum Gravity And Generalized Matrix Models,''
  Mod.\ Phys.\ Lett.\  A {\bf 6}, 1133 (1991).
  \doi{10.1142/S0217732391001184}
\\
  M.~Gross,
  ``Tensor models and simplicial quantum gravity in $>$ 2-D,''
  Nucl.\ Phys.\ Proc.\ Suppl.\  {\bf 25A}, 144 (1992).
  \doi{10.1016/S0920-5632(05)80015-5}
\\
  N.~Sasakura,
  ``Tensor model for gravity and orientability of manifold,''
  Mod.\ Phys.\ Lett.\  A {\bf 6}, 2613 (1991).
  \doi{10.1142/S0217732391003055}

\bibitem{GurauBook}
  R.~Gurau,
  ``Random tensors,''
  Oxford University Press 2016.

\bibitem{ItalianSurvey}
  M.~Ferri, C.~Gagliardi and L.~Grasselli,
  ``A graph-theoretical representation of PL-manifolds -- A survey on crystallizations,''
  Aequationes Mathematicae 31 (1986) 121.
  \doi{10.1007/BF02188181}

\bibitem{LinsMandel}
  S.~Lins and A.~Mandel,
  ``Graph-encoded 3-manifolds,''
  Discrete Mathematics 57 (1985) 261.
  \doi{10.1016/0012-365X(85)90179-7}

\bibitem{Balanced}
  L.~J.~Billera and A.~Björner,
  ``Face numbers of polytopes and complexes,''
  Handbook of discrete and computational geometry (Jacob E. Goodman and Joseph O?Rourke, eds.), CRC Press Series on Discrete Mathematics and its Applications, CRC Press, 1997, pp. 407?430\\
  I.~Izmestiev and M.~Joswig,
  ``Branched coverings, triangulations, and 3-manifolds,''
  Adv. Geom. {\bf 3}(2):191?225, 2003
  \doi{10.1515/advg.2003.013},
  \arxiv{math/0108202}\\
  M.~Joswig,
  ``Projectivities in simplicial complexes and colorings of simple polytopes,''
  Math. Z. {\bf 240}(2):243?259, 2002
  \doi{10.1007/s002090100381},
  \arxiv{math/0102186}\\
  S.~Klee and I.~Novik,
  ``Lower Bound Theorems and a Generalized Lower Bound Conjecture for balanced simplicial complexes,''
  Mathematika {\bf 62} (2016) 441?477,
  \doi{10.1112/S0025579315000297},
  \arxiv{1409.5094}
  R.~P.~Stanley, 
  ``Combinatorics and commutative algebra,''
  Progress in Mathematics, vol. 41, Birkhaüser Boston Inc., Boston, MA, 1983. MR MR725505 (85b:05002)
  


\bibitem{Uncoloring}
  V.~Bonzom, R.~Gurau and V.~Rivasseau,
  ``Random tensor models in the large N limit: Uncoloring the colored tensor models,''
  Phys.\ Rev.\ D {\bf 85}, 084037 (2012)
  \doi{10.1103/PhysRevD.85.084037},
  \arxiv{1202.3637} [hep-th].

\bibitem{Regge} 
  T.~Regge,
  ``General Relativity Without Coordinates,''
  Nuovo Cim.\  {\bf 19}, 558 (1961).
  \doi{10.1007/BF02733251}

\bibitem{GurauSchaeffer}
  R.~Gurau and G.~Schaeffer,
  ``Regular colored graphs of positive degree,''
  \arxiv{1307.5279} [math.CO].

\bibitem{1/NExpansion}
  R.~Gurau,
  ``The complete 1/N expansion of colored tensor models in arbitrary dimension,''
  Annales Henri Poincar\'e {\bf 13}, 399 (2012)
  \doi{10.1007/s00023-011-0118-z},
  \arxiv{1102.5759} [gr-qc].

\bibitem{Melons}
  V.~Bonzom, R.~Gurau, A.~Riello and V.~Rivasseau,
  ``Critical behavior of colored tensor models in the large N limit,''
  Nucl.\ Phys.\ B {\bf 853} (2011) 174
  \doi{10.1016/j.nuclphysb.2011.07.022},
  \arxiv{1105.3122} [hep-th].

\bibitem{MOClassification}
  E.~Fusy and A.~Tanasa,
  ``Asymptotic expansion of the multi-orientable random tensor model,''
   Electron. J. Combin. 22 (1) (2015), \#P1.52
  \arxiv{1408.5725} [math.CO].

\bibitem{Meanders}
  V.~Bonzom and F.~Combes,
  ``The calculation of expectation values in Gaussian random tensor theory via meanders,''
  AIHP-D {\bf 1} (2014) 443-485,
  \doi{10.4171/AIHPD/13},
  \arxiv{1310.3606}.

\bibitem{MelonoPlanar} 
  V.~Bonzom, T.~Delepouve and V.~Rivasseau,
  ``Enhancing non-melonic triangulations: A tensor model mixing melonic and planar maps,''
  Nucl.\ Phys.\ B {\bf 895}, 161 (2015)
  \doi{10.1016/j.nuclphysb.2015.04.004},
  \arxiv{1502.01365} [math-ph].

\bibitem{O(N)Model}
  S.~Carrozza and A.~Tanasa,
  ``$O(N)$ Random Tensor Models,''
  Lett.\ Math.\ Phys.\  {\bf 106}, no. 11, 1531 (2016)
  \doi{10.1007/s11005-016-0879-x}
  \arxiv{1512.06718}.
  
\bibitem{StuffedWalshMaps}
  V.~Bonzom, L.~Lionni and V.~Rivasseau,
  ``Colored triangulations of arbitrary dimensions are stuffed Walsh maps,''
  \arxiv{1508.03805} [math.CO].

\bibitem{Catalogues}
  M.~R.~Casali and C.~Gagliardi,
  ``A code for m-bipartite edge-coloured graphs,''
  Rend. Ist. Mat. Univ. Trieste 32 suppl.1, (2001), pp. 55?76.\\
  A.~Marani, M.~Rivi, P.~Cristofori, 
  ``Generation of Catalogues of PL n-manifolds: Computational Aspects on HPC Systems,''
  14-th International Symposium on Symbolic and Numeric Algorithms for Scientific Computing, Timisoara (2012).

\bibitem{FusyThese}
  E.~Fusy,
  ``Combinatorics of planar maps and algorithmic applications,''
  PhD manuscript
  \url{http://www.lix.polytechnique.fr/Labo/Eric.Fusy/Theses/these_eric_fusy.pdf}
  
  
\bibitem{DartoisDoubleScaling} 
  S.~Dartois, R.~Gurau and V.~Rivasseau,
  ``Double Scaling in Tensor Models with a Quartic Interaction,''
  JHEP {\bf 1309}, 088 (2013)
  \doi{10.1007/JHEP09(2013)088},
  \arxiv{1307.5281} [hep-th].

\bibitem{Universality}
  R.~Gurau,
  ``Universality for Random Tensors,''
  Ann.\ Inst.\ H.\ Poincare Probab.\ Statist.\  {\bf 50}, no. 4, 1474 (2014)
  \doi{10.1214/13-AIHP567},
  \arxiv{1111.0519} [math.PR].

\end{thebibliography}
\end{document}